\documentclass[]{article}
\usepackage{amsmath,amsfonts,amsthm,url,graphicx,enumerate}
\usepackage{xcolor}
\usepackage{subfigure}
\newtheorem{exercise}{Exercise}
\newtheorem{theorem}{Theorem}
\newtheorem{fact}{Fact}
\newtheorem{remark}{Remark}
\newtheorem{solution}{Solution}
\newcommand{\RR}{\mathbb{R}}
\newcommand{\Tt}{\mathbb{T}_T}
\newcommand{\bt}{\bar{t}}
\newcommand{\melnikovconds}{ {\em M.1-M.2}}
\begin{document}
\title{Notes on the computation of periodic orbits using Newton and
Melnikov's method}
\date{February 2017}
\author{Albert Granados\\
algr@dtu.dk\\
Department of Applied Mathematics and
Computer Science,\\ Technical University of Denmark,\\
Building 303B, 2800
Kgns. Lyngby, Denmark}
\maketitle
\tableofcontents
\section{Introduction}
The main goal of these sessions will be to understand the dynamics of a
system exhibiting a continuum of periodic orbits when we add a small
periodic forcing. The most paradigmatic example is probably the
perturbed pendulum; however, such systems massively appear in real
applications, specially in celestial mechanics. This has given rise to
classical problems exhibiting extremely rich dynamics, such as the
restricted three body problem.

In these sessions we will see some theoretical results, but we will
mainly visualize them through analytical and numerical computations
for a particular example, the forced pendulum:
\begin{equation}
\ddot{u}+\sin(u)=\varepsilon g(x,t),
\label{eq:perturbed_pendulum}
\end{equation}
where $\varepsilon\ge0$ is a small parameter, $x=(u,\dot{u})$ and
$g(x,t)$ a periodic forcing: $g(x,t+T)=g(x,t)$.
\begin{figure}
\begin{center}
\includegraphics[angle=-90,width=0.7\textwidth]{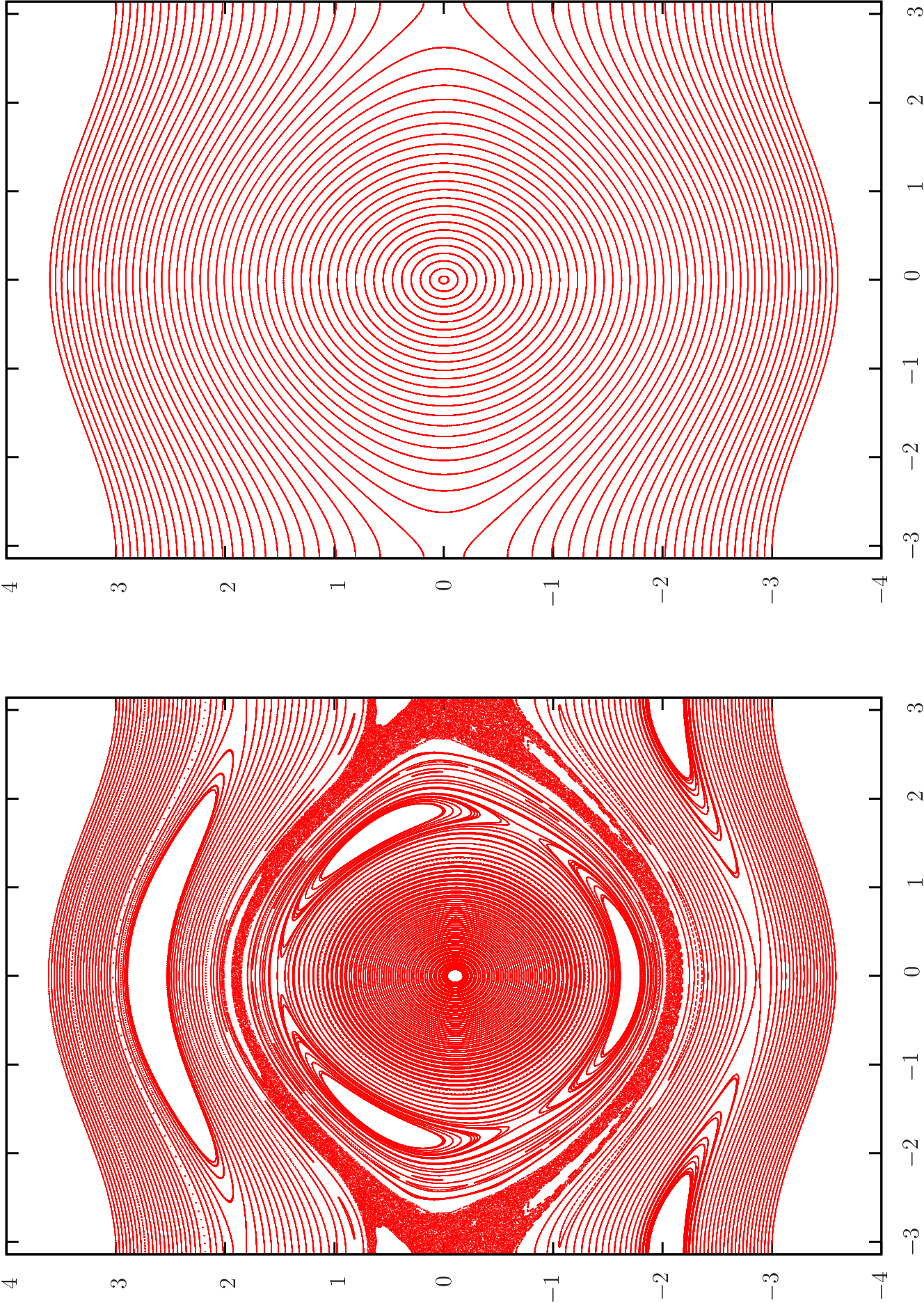}
\end{center}
\caption{Phase portrait of the unperturbed pendulum (right). Phase
portrait for $\varepsilon>0$ (left).}
\label{fig:pendulum}
\end{figure}
We will mainly see that, when $\varepsilon=0$ the face portrait looks
like Figure~\ref{fig:pendulum} (right) and that, when $\varepsilon>0$
looks like Figure~\ref{fig:pendulum} (left). We will learn theoretical
and numerical techniques to compute the surviving resonant periodic
orbits (large ``holes'' in Figure~\ref{fig:pendulum} left).\\
Although everything we will see will be generic for any $g(x,t)$, we
will fix this forcing from now on say to
\begin{equation*}
g(x,t)=\sin(\omega t),
\end{equation*}
ant $T$ becomes $T=\frac{2\pi}{\omega}$.\\
It will be useful to write the system as a first order one by
increasing the dimension. Let $v=\dot{u}$, then we can write it as
\begin{equation}
\left\{
\begin{aligned}
\dot{u}&=v\\
\dot{v}&=-\sin(u)+\varepsilon \sin(\omega t).
\end{aligned}
\right.
\label{eq:fo_system}
\end{equation}
\section{Dynamics of the pendulum around the elliptic
point}\label{sec:quasi_integrable_dynamics}
\subsection{Unperturbed phase portrait}
We will start with some exercises in order to understand the dynamics
of the unperturbed pendulum ($\varepsilon=0$).
\begin{exercise}
Consider the unperturbed equations of the pendulum by setting
$\varepsilon=0$ in Eq.~\eqref{eq:fo_system}:
\begin{equation}
\left\{
\begin{aligned}
\dot{u}&=v\\
\dot{v}&=-\sin(u)
\end{aligned}
\right.
\label{eq:ode_pendulum}
\end{equation}
\begin{enumerate}
\item Obtain the equilibrium points and state their type.
\item Obtain a Hamiltonian for system~\eqref{eq:ode_pendulum}.
\end{enumerate}
\end{exercise}

\begin{solution}
\begin{enumerate}
\item The system possesses an elliptic equilibrium point at the
origin and two saddle points at $( \pm\pi,0)$. 
\item The function
\begin{equation*}
H(u,v)=\frac{v^2}{2}-\cos(u)
\end{equation*}
is the Hamiltonian of system~\eqref{eq:ode_pendulum} up to a constant value. 
\end{enumerate}
\end{solution}

\begin{exercise}
Write a small program in Matlab to plot the phase portrait of
system~\eqref{eq:ode_pendulum}.
\end{exercise}
\begin{solution}
The phase portrait for the unperturbed pendulum looks like the one in
Figure~\ref{fig:pendulum} (right). The phase portrait is divided in
two regions by two separatrices. These separatrices are given by two
homoclinic trajectories connecting the two saddle points, and are
located at the level of energy $H(u,v)=\frac{0^2}{2}-\cos(\pi)=1$. In
the inner part we find the elliptic equilibrium point, which is
surrounded by a continuum of periodic orbits. Provided that the energy
level of the elliptic point is $H(0,0)=-1$, each periodic orbit is
located a certain level of energy between $1$ and $-1$: $H(u,v)=c$,
$c\in(-1,1)$.  Outside the homoclinic loop we also find a continuum of
periodic orbits, given by $H(u,v)=c$, $c>1$.
\end{solution}
From now on we will focus on what happens with the periodic orbits
inside the homoclinic loop when we add a small periodic forcing to
System~\eqref{eq:ode_pendulum} (say $\sin(\omega t)$). To this end, it
will be crucial to compute the periods of the periodic orbits for
$\varepsilon=0$.  Let's us now see a general theoretical approach to
obtain a formula.

Assume we have Hamiltonian system of the form
\begin{equation*}
H(u,v)=\frac{v^2}{2}+V(u),
\end{equation*}
where $V(u)$ is a potential energy. Assume that we have a periodic
orbit at level of energy $c\in \RR$ given by $H(u,v)=c$.
Then, we can compute its period, $T_c$, as follows:
\begin{align*}
T_c&=\oint_{H(u,v)=c}dt=\oint_{H(u,v)=c}\frac{1}{\frac{du}{dt}}du\\
&=\oint_{H(u,v)=c}\frac{1}{\dot{u}}du.
\end{align*}
using that the system is Hamiltonian and, hence,
$\dot{u}=\frac{\partial H(u,v)}{\partial v}=v$ we get
\begin{equation*}
T_c=\oint_{H(u,v)=c}\frac{1}{v}du.
\end{equation*}
Isolating $v$ from $H(u,v)=c$,
\begin{equation*}
T_c=\oint_{H(u,v)=c}\frac{1}{\sqrt{2(c-V(u))}}du.
\end{equation*}

But now we have two problems. First, solving such an integral
explicitly can be a nightmare, if possible at all. Second, doing
integrals numerically is a difficult task, it's slow and imprecise.
Moreover, note that, although the integral converges, the integrand
has an asymptote at the crossing point with the horizontal axis:
$(0,u_1)$ such that $V(u_1)=c$. This is specially problematic when
computing this integral numerically.\\
Fortunately, there is an alternative: compute a Poincar\'e map using a
transversal section to the periodic orbit and capture the return time.

\begin{figure}
\begin{center}
\includegraphics[width=1\textwidth]{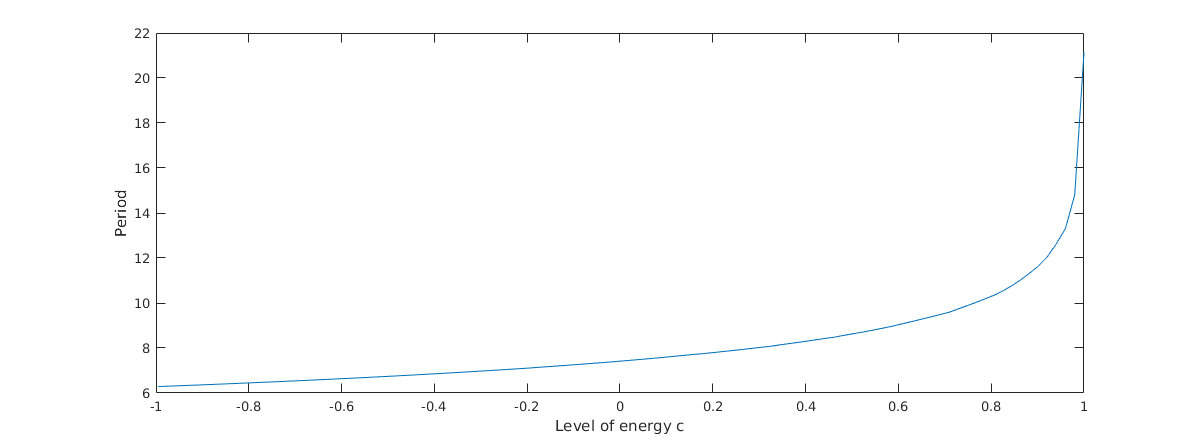}
\end{center}
\caption{Periods of the periodic orbits for the unperturbed pendulum.}
\label{fig:unper_periods}
\end{figure}

\begin{exercise}\label{ex:period}
Write a small program in Matlab to compute the periods of the periodic
orbits using the Poincar\'e map from the section $\left\{ u=0 \right\}$
to itself. 
\end{exercise}
\begin{solution}
Use the function {\em period.m} to obtain something similar to
Figure~\ref{fig:unper_periods}. Note that, the argument of this
function is $v_0$ and starts integrating at $(0,v_0)$. Hence, the
corresponding level of energy becomes $H(0,v_0)=\frac{v_0^2}{2}-1$.\\
Note that there is a vertical asymptote at $c=1$, which corresponds to
the energy level of the homoclinic loop.
\end{solution}

\subsection{Stroboscopic vs Poincar\'e map}
In general, for autonomous systems, in order to study the existence of
periodic orbits one typically considers a Poincar\'e map:
\begin{equation*}
P: \Sigma:\longrightarrow \Sigma,
\end{equation*}
where $\Sigma$ is a co-dimension one section transversal to the
periodic orbits. In our case, we could take the vertical axis
$\Sigma=\left\{ (u,v)\,|\,u=0 \right\}$. Then, fixed points of the
Poincar\'e map, points $x_0\in \Sigma$ such that $P(x_0)=x_0$, become
initial conditions for periodic orbits for the flow, whose period is
the flying time. Periodic points of the Poincar\'e map, $P^n(x_0)=x_0$,
also give rise to periodic orbits for the flow, which cross $n$ times
the Poincar\'e section and their period becomes the addition of all
flying times between consecutive impacts with the Poincar\'e section
before $x_0$ is reached again.\\
However, if the system is not autonomous (but periodic), then one
needs to take into account the initial time and consider Poincar\'e
sections of the form $\Sigma\times \Tt$, with $\Tt=\RR/T\mathbb{Z}$.
Then, a sufficient condition for the existence of a periodic orbits
becomes $P^n(x_0,t_0)=(x_0,t_0+mT)$. That is, after the point $x_0$ is
reached again after $n$ crossings and the total spent time is a
multiple of $T$. Then, the flow possesses an $mT$-periodic orbit
crossing the section $n$ times. Note that, if total time spent to
reach $x_0$ again is not a multiple of $T$, then nothing can be said
about the existence of a periodic orbit.

One of inconveniences of using Poincar\'e maps with non-autonomous
systems is the need of computing the flying time. Although this can be
done numerically, it requires extra computations than simply
numerically integrating a flow, as one needs to compute the crossing
with the section. Alternatively, one can use the stroboscopic map,
which consists of integrating the system for a time $T$:
\begin{equation}
\begin{array}{cccc}
s_\varepsilon:&\Sigma_{t_0}&\longrightarrow &\Sigma_{t_0}\\
&x&\longmapsto&\varphi_\varepsilon(T;x,t_0),
\end{array}
\label{eq:strobomap}
\end{equation}
where
\begin{equation*}
\Sigma_{t_0}=\left\{ (x,t)\in \RR^2\times \Tt
\right\},\,\Tt:=\RR/T\mathbb{Z}.
\end{equation*}
Then, if $s_\varepsilon^m(x_0)=x_0$, then $x_0$ is the initial condition for a
$mT$-periodic orbit:
\begin{equation*}
s_\varepsilon^m(x_0)=\varphi_\varepsilon(t_0+mT;x_0,t_0)=x_0.
\end{equation*}

\begin{exercise}\label{ex:strobo_map}
Use the Matlab script {\em iter\_strobo.m} or {\em
iter\_strobo\_parallel.m} to simulate the stroboscopic map for a set
of different initial conditions. The latter parallelizes the loop and
is faster if you have more cores available.
\end{exercise}
\subsection{Dynamics of the stroboscopic map: unperturbed case}\label{sec:unpert_strobo}
As mentioned above, for $\varepsilon=0$, System~\eqref{eq:fo_system} possesses a
continuum of periodic orbits surrounding the origin and whose periods
increase when approaching the homoclinic loop. Therefore, when picking
up a point, all its iterates will be located in an invariant curve,
which is nothing else than the periodic orbit of the time-continuous
system. However, the dynamics of the map will exhibit different
properties depending on the relation between $T$ and $T_c$.\\
Imagine we iterate a point $x_0=(u_0,v_0)$ such that $H(x_0)=c\in
(-1,1)$. Then, if $T$ and $T_c$ are congruent ($T$ and $T_c$ are
rational multiples), this point will be a periodic point of the
stroboscopic map. Let us assume that there exists some
$m\in\mathbb{Z}$ such that $T_c=mT$. Then, the stroboscopic map will
return to $x_0$ after $m$ iterates:
\begin{equation*}
s_0^m(x_0)=x_0.
\end{equation*}
Indeed, all points at the level of energy $c$ will be $m$-periodic points
\begin{figure}
\begin{center}
\includegraphics[width=0.8\textwidth]{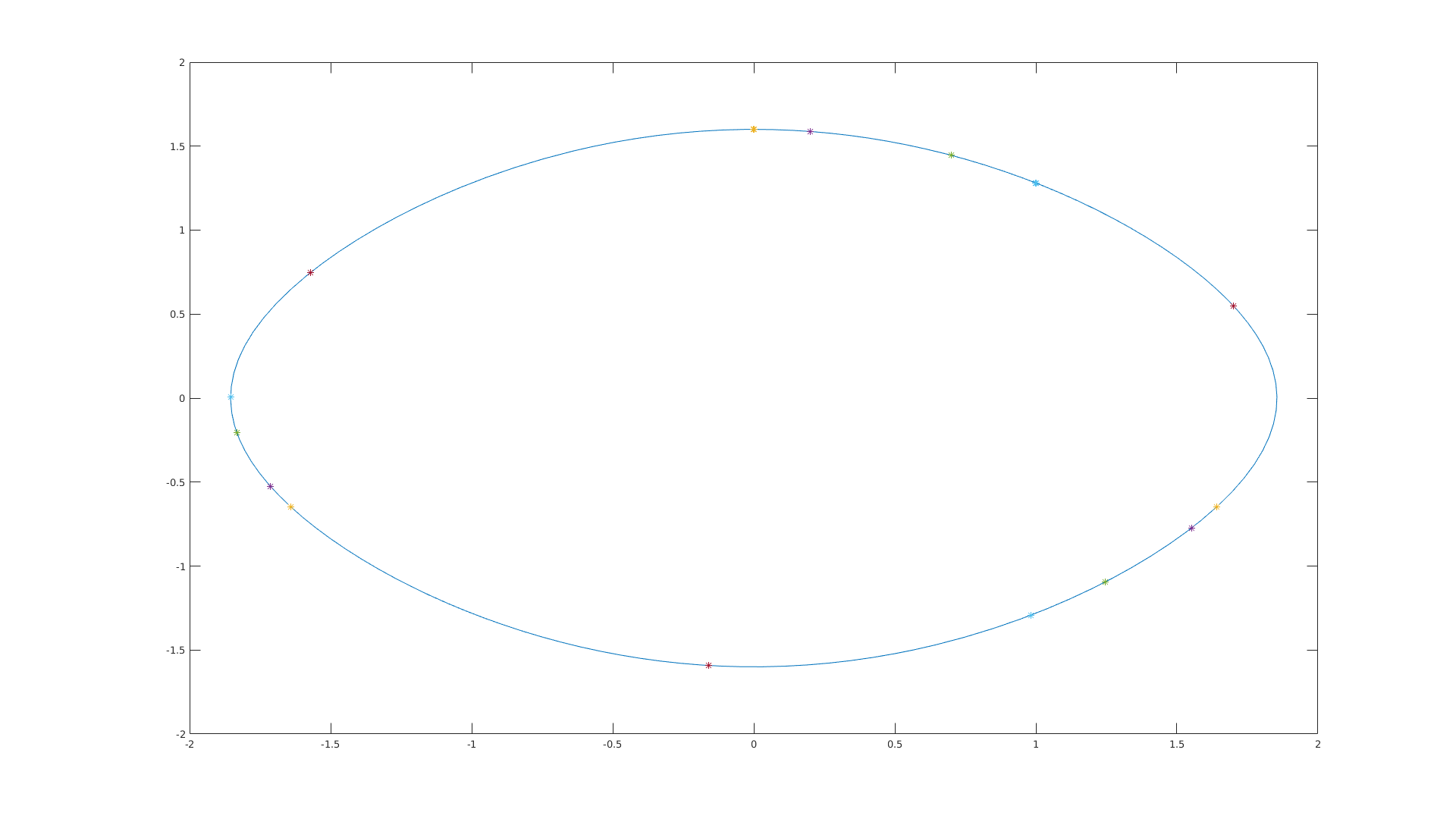}
\end{center}
\caption{Continuous curve: periodic orbit with period $T_c=3T$.
Points: different initial conditions (with different colors) and its iterates.}
\label{fig:unpert3po}
\end{figure}
of the stroboscopic map. This is illustrated in
Figure~\ref{fig:unpert3po} we show a periodic orbit with $T_c=3T$ and
the iterates of different initial conditions taken at this trajectory.
Note that they all form a $3$-periodic orbit: for any $x_0$ at that
level of energy we have $s_0^3(x_0)=x_0$.  The parameters have been
bound as follows. We have chosen an arbitrary level of energy
$H(0,1.6)=c$. Then, using the function {\em period.m}, we have found
the period of this level set, $T_c$, and then chosen $T$ such that
\begin{equation*}
T_c=3T.
\end{equation*}
Note that, this relation imposes the fact that, after iterating
$3$ times the stroboscopic map, the time-continuous system has
performed exactly one loop. However, let us see what happens if we
choose $T$ such that $T_c=\frac{mT}{n}$.
\begin{figure}
\begin{center}
\includegraphics[width=0.8\textwidth]{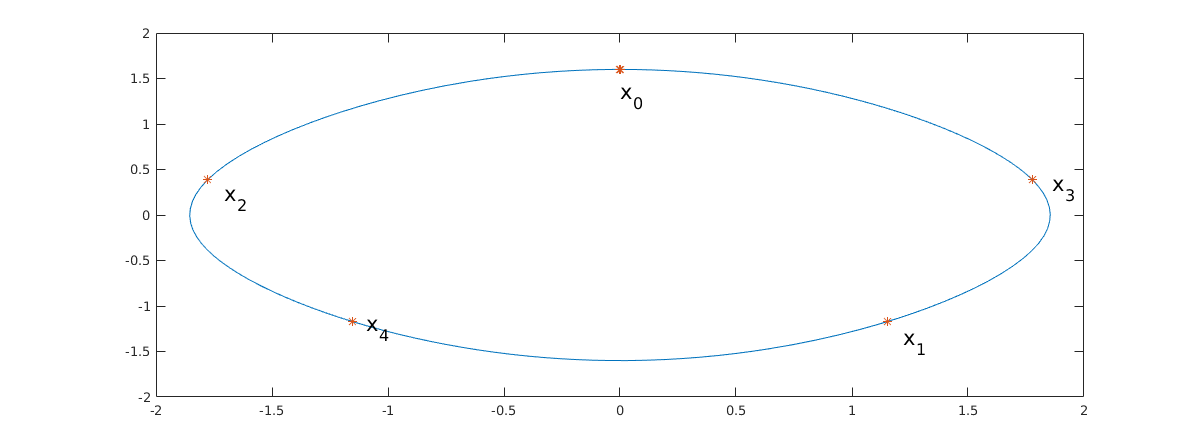}
\end{center}
\caption{Continuous curve: periodic orbit with period
$T_c=\frac{2T}{5}$.
Red points: iterates of the point $(0,v_0)$ such that $H(0,v_0)=c$ and
$T_c=\frac{2T}{5}$.}
\label{fig:unpert2-5po}
\end{figure}
In Figure~\ref{fig:unpert2-5po} we show the situation at the same
level of energy when choosing $T$ such that $T_c=\frac{2T}{5}$. In
this case we only show the iterates of a unique point at that level
energy, but all points have the same behaviour. Note that, in this
case, the points are not ordered with increasing angle, but at each
iteration one skips one point ($s(x_i)=x_{i+1}$). Therefore, when the
loop is closed, we have twisted twice along the periodic orbit of the
time continuous system. This is called a subharmonic orbit and one
says that has rotation number $\frac{2}{5}$. Roughly speaking, the
rotation number of a periodic orbit measures the average angle twisted
at each iterate. Note that, in the previous example, one has a
harmonic periodic orbit with rotation number $\frac{1}{3}$.

Summarizing we have the following fact.
\begin{fact}
If $T$ is such that $T_c=\frac{mT}{n}$, then all points satisfying
$H(x_0)=c$ are $m$ periodic points of the stroboscopic map:
\begin{equation*}
s_0^m(x_0)=x_0.
\end{equation*}
When the cycle is closed, the iterates have performed $n$ loops around
the origin. In other words, the rotation number of the periodic orbit
is $n/m$.
\end{fact}
Let us now discuss the situation when $T$ and $T_c$ are
incommensurable (such $m$ and $n$ do not exist). Recalling that
the time-continuous system is $T$-periodic, any periodic orbit needs
to spend a time multiple of $T$ to close. If there
doesn't exist any integer combination such that $nT_c=mT$, then any
point at the level of energy $H(u,v)=c$ cannot be closed by iterating
the stroboscopic map. To see that, let us assume that $s_0^m(x_0)=x_0$
for some $m$. Then there is a point satisfying $H(x_0)=c$ such that
\begin{equation*}
\varphi_0(mT;x_0)=x_0.
\end{equation*}
On the other hand, we also have that $\varphi_0(T_c;x_0)=x_0$ and,
therefore, $T_c=mT$, which is a contradiction, and the point never
returns to itself when iterated by the stroboscopic map.\\
As a consequence, as the periodic orbit is an invariant curve for the
stroboscopic map, when iterating the map we will densely fill this
curve. In theory, if the map restricted to the invariant curve is not
$C^2$, this can lead to a Denjoy counterexample (see~\cite{Nit71}) and one could get a
Cantor set instead of filling this curve densely.
\begin{figure}
\begin{center}
\includegraphics[width=0.8\textwidth]{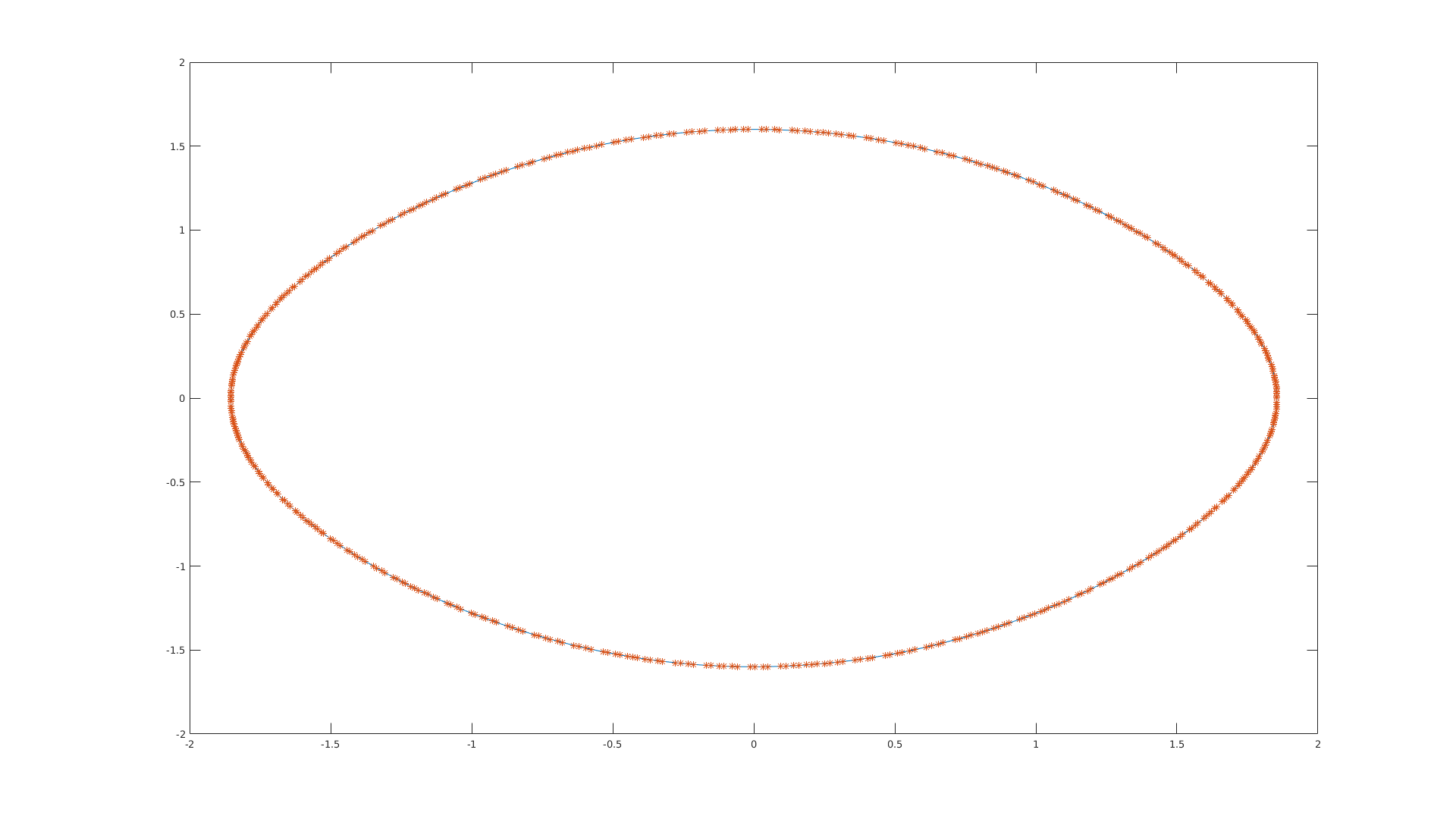}
\end{center}
\caption{Blue curve: periodic orbit with period
$T_c=\frac{T}{\sqrt{2}}$.  Red points: $500$ iterates of the point
$(0,v_0)$ such that $H(0,v_0)=c$ and $T_c=\frac{T}{\sqrt{2}}$.}
\label{fig:unpertsqrt2}
\end{figure}
In Figure we show the iterates of a point belonging at the level of
energy $H(0,1.6)=c$ with $T_c=\frac{T}{\sqrt{2}}$.
\begin{remark}
For $\varepsilon=0$ the system is autonomous. Therefore, the
dynamics of $s(x_0)$ is independent on the election of $t_0$ in
Equation~\eqref{eq:strobomap}. 
\end{remark}
\begin{exercise}
Use the code written in Exercise~\ref{ex:strobo_map} to simulate the
stroboscopic map for $\varepsilon=0$ in different situations. Fix
$t_0=0$ and some level of energy $H(u,v)=c$. Compute $T_c$ using the
function developed in Exercise~\ref{ex:period} and then
\begin{enumerate}[i)]
\item choose $T$ so that it is congruent with $T_c$. Try different
combinations as described before to obtain equivalent plots as in
Figures~\ref{fig:unpert3po} and~\ref{fig:unpert2-5po}.
\item choose $T$ so that it is unconmeasurable with $T_c$ and plot
$500$ iterates of a point in that level of energy.
\item Repeat the previous experiments with another value of $t_0$.
What differences do you observe?
\end{enumerate}
\end{exercise}
\begin{solution}
What you should observe is explained in the text. Note that you will
need to mofify {\em period.m} to include the initial integrating time
$t_0$.
\end{solution}
\subsection{Dynamics of the stroboscopic map: perturbed case}\label{sec:pert_strobomap}
When activating $\varepsilon>0$ (small), the dynamics of
System~\eqref{eq:fo_system} become very rich. Typically one can observe homoclinic
tangles, horseshoes leading to chaos, invariant curves and periodic
orbits. We will focus on the existence of periodic orbits and, later
on, on its computation.\\
\unitlength=\textwidth
\begin{figure}
\begin{center}
\includegraphics[width=1\textwidth]{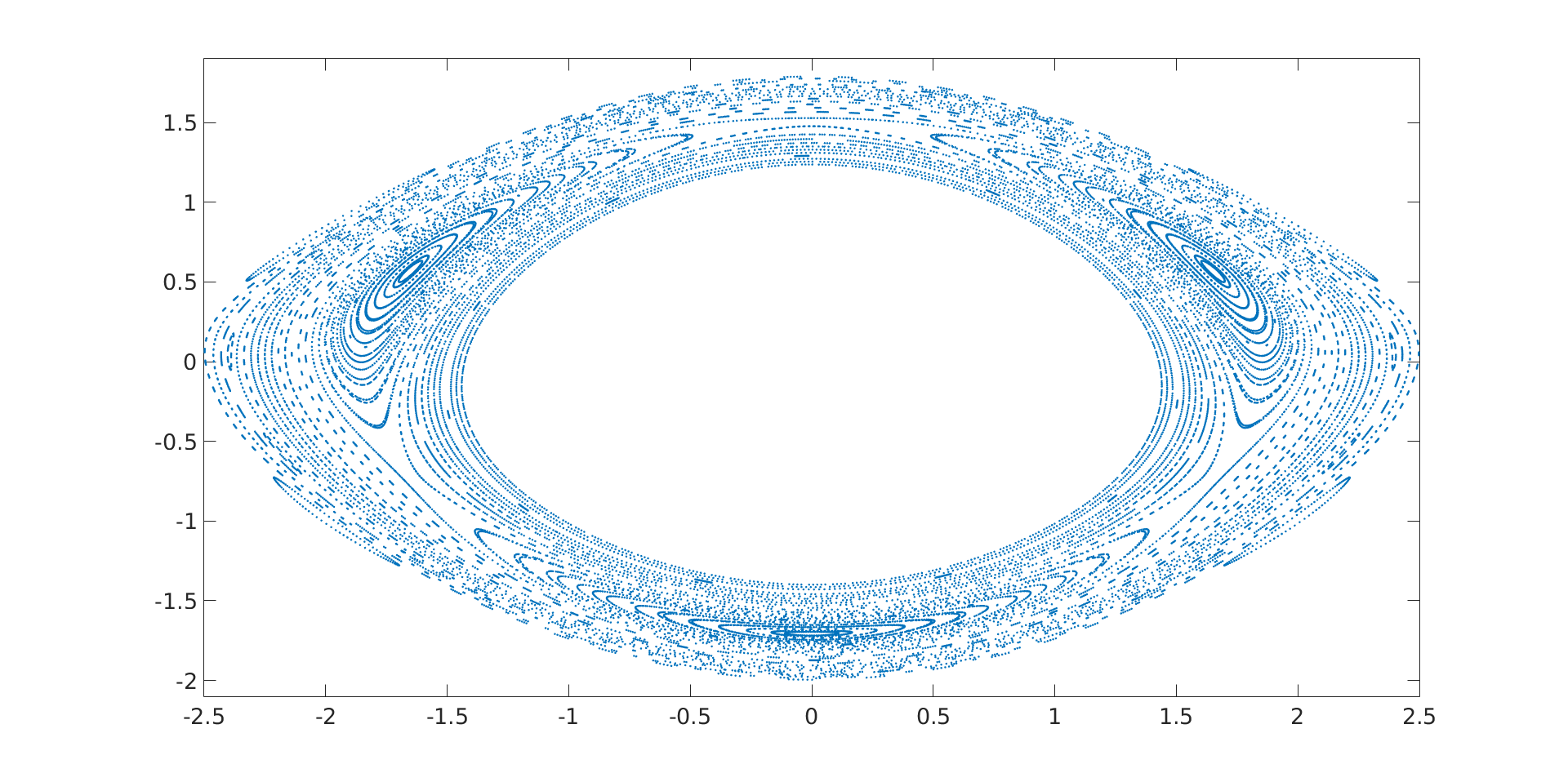}
\end{center}
\caption{Dynamics of the stroboscopic map for $\varepsilon=0.2$ close
to a $3$-periodic orbit for $t_0=0$. The figure has been obtained by iterating
initial conditions along a line going through the elliptic periodic
orbit surrounding the level of energy
$H(0,1.6)$ and $T=T_c/3$.}
\label{fig:pert3pot00}
\end{figure}
In Figure~\ref{fig:pert3pot00} we show the dynamics of the stroboscopic
map for $\varepsilon=0.2$ close to the level of energy $c=H(0,1.6)$
(chosen arbitrarily in Section~\ref{sec:unpert_strobo}) and $t_0=0$. When choosing
$T=T_c/3$, one observes two $3$-periodic orbits. One periodic orbit is
of the saddle type, and one iterate is close to the vertical axis. The
other periodic orbit is of the elliptic type, as the three points of
this periodic orbit are surrounded by (apparent) curves.\\
\begin{exercise}\label{ex:perturbed_strobo_map}
Use the code written in Exercise~\ref{ex:strobo_map} to simulate the
stroboscopic map of System~\eqref{eq:fo_system} for $\varepsilon>0$
(small) in the following situations.
\begin{enumerate}
\item set $t_0$ and iterate initial conditions around a resonance of the
unperturbed system. That is, choose $T=T_c/m$ for
some $c\in(-1,1)$, iterate use some initial conditions along
$(0,v_0-2\varepsilon)$ to $(0,v_0+2\varepsilon)$, with $c=\frac{v_0^2}{2}-1$. Alternatively,
you can take initial conditions along a line (the diagonal for
example): $(v_1-2\varepsilon,v_1-2\varepsilon)$ to
$(v_1+2\varepsilon,v_1+2\varepsilon)$, with
$H(v_1,v_1)=\frac{v_1^2}{2}-\cos{v_1}$.  How many periodic orbits do you observe?
\item Repeat the same simulation choosing some $t_0=0$. What changes
do you observe? What if you choose $t_0=T$?
\end{enumerate} 
\end{exercise}
\begin{figure}
\begin{center}
\includegraphics[width=1\textwidth]{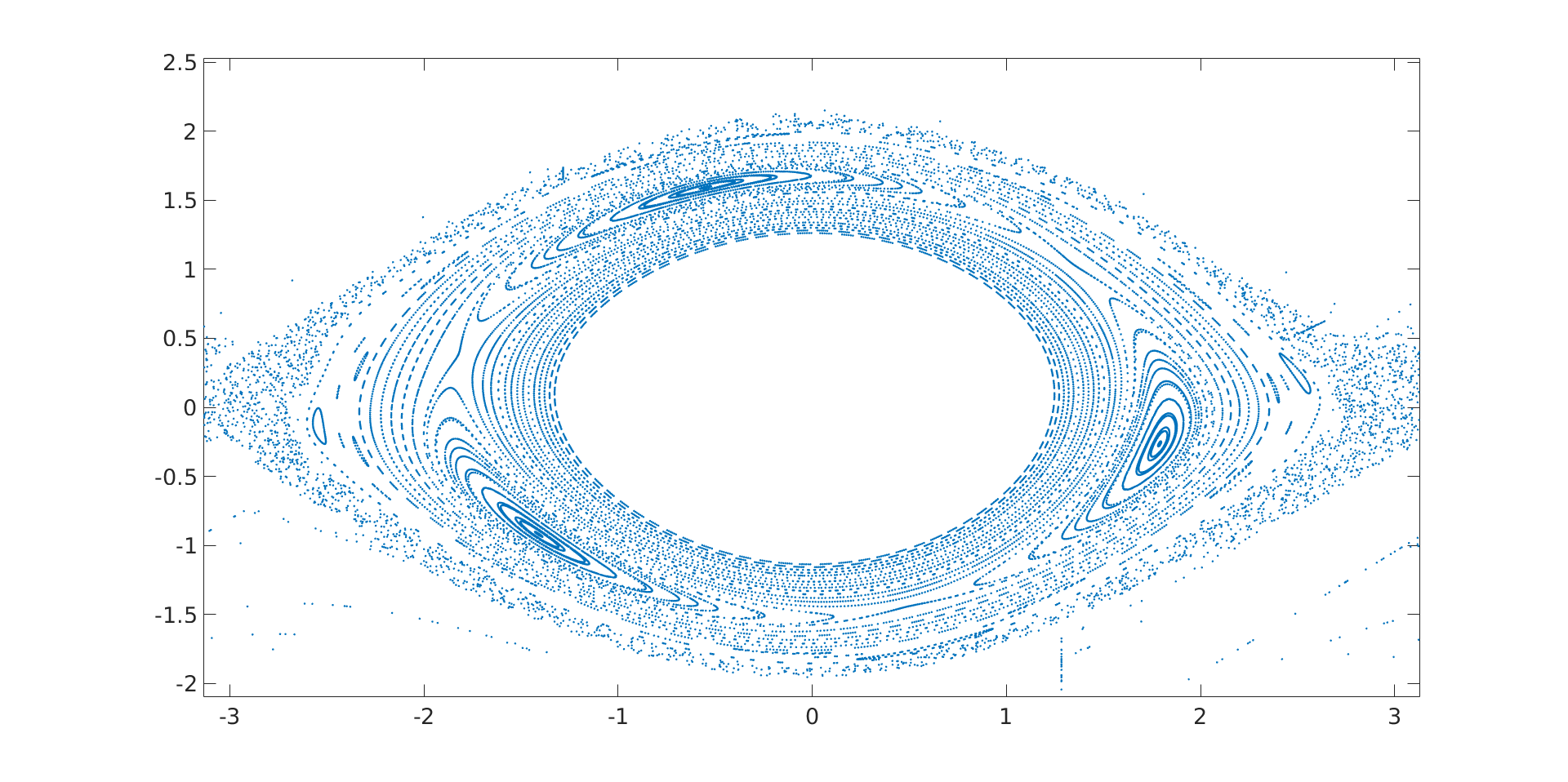}
\end{center}
\caption{Dynamics of the stroboscopic map for $\varepsilon=0.2$ close
to a $3$-periodic orbit for $t_0=1$. The figure has been obtained by iterating
initial conditions along a line going through the elliptic periodic
orbit surrounding the level of energy
$H(0,1.6)$ and $T=T_c/3$.}
\label{fig:pert3pot01}
\end{figure}
\begin{solution}
In Figure~\ref{fig:pert3pot00} we have what you should observe for
$v_0=1.6$, $T=T_c/3$, $t_0=0$ and taking initial conditions along the
line $y=0.5/1.6x$. In Figure~\ref{fig:pert3pot01} we have chosen
$t_0=1$ and initial conditions along a straight line going through the
elliptic periodic orbit. Note that the dynamics is the same as before
but the portrait is twisted. As the system is $T$-periodic, one should
observe the same for $t_0=T$ than for $t_0=0$. 
\end{solution}

\begin{remark}
As for $\varepsilon>0$ System~\eqref{eq:fo_system} is non-autonomous,
the dynamics of the stroboscopic map depend on the election of
$t_0$. However, the dynamics are qualitatively the same but twist when
$t_0$ is varied. Compare Figures~\ref{fig:pert3pot00}
and~\ref{fig:pert3pot01}.
\end{remark}

\section{The Melnikov method for subharmonic periodic orbits}
In Section~\ref{sec:pert_strobomap} we have seen using numerical
simulation that, when $T_c$ and $T$ are congruent, the corresponding
periodic orbit may persist for $\varepsilon>0$. In this section we
present a result (the Melnikov method) that provides sufficient conditions for the
persistence of periodic orbits resonant with the period of the
perturbation. Indeed, it does not only predict this persistence but
also provides information regarding the location of these periodic
orbits. As we will see in Section~\ref{sec:compute_po_strobo}, this
will become very useful when computing periodic orbits using Newton's
method.
\subsection{Persistence of periodic orbits}
The Melnikov method can be applied to a larger class of systems than
the pendulum. More precisely, we consider a general planar field of
the form
\begin{equation}
\dot{x}=f(x)+\varepsilon g(x,t)=
\left(
\begin{array}[]{c}
f_1(x)\\f_2(x)
\end{array}
 \right)
+\varepsilon\left( 
\begin{array}[]{c}
g_1(x,t)\\g_2(x,t)
\end{array}
 \right)
,
\label{eq:field_melnikov}
\end{equation}
where $x=(u,v)\in\RR^2$, $\varepsilon\ge0$ is a small parameter and
$g(x,t)$ is $T$-periodic in $t$:
\begin{equation*}
g(x,t+T)=g(x,t), \forall x\in \RR^2.
\end{equation*}
For simplicity, let us assume that, for $\varepsilon=0$, the unperturbed system
\begin{equation*}
\dot{x}=f(x)
\end{equation*}
is Hamiltonian. That is, there exists a function $H(u,v)$ such that
\begin{align*}
\dot{u}&=\frac{\partial H(u,v)}{\partial v}\\
\dot{v}&=-\frac{\partial H(u,v)}{\partial u}.
\end{align*}
Moreover, let us assume  that for $\varepsilon=0$, the system
satisfies the following conditions:
\begin{enumerate}[{\em M.1}]
\item There exists a compact region completely covered by a continuum of
periodic orbits.
\item Let $T_c$ be the period of the periodic orbit located at
the energy level $c$:
\begin{equation*}
\left\{ (u,v)\in\RR^2,\,H(u,v)=c \right\}.
\end{equation*}
Assume that $\frac{d}{dc}T_c\neq 0$.
\end{enumerate}
Then we have the following result (\cite{GucHol83,Mel63}):
\begin{theorem}\label{theo:Melnikov}
Assume the above conditions are satisfied. Assume that the
unperturbed system has a periodic orbit, $H(u,v)=c$, of period
\begin{equation}
T_c=\frac{m}{n}T,
\label{eq:congruency}
\end{equation}
where $T$ is the period of the periodic forcing. Let $x_0=(u_0,v_0)$
be an initial condition for such a periodic orbit (that is,
$H(x_0)=c$), and let 
\begin{equation*}
\varphi_0(t;x_0)
\end{equation*}
be the unperturbed flow at such initial condition. Then let us define
the so-called (subharmonic) Melnikov function
\begin{equation}
\begin{aligned}
M(t_0)&=\int_{0}^{mT}f(\varphi_0(t;x_0))\wedge
g(\varphi_0(t;x_0),t+t_0)dt\\
&=\int_{0}^{mT}\Big(f_1(\varphi_0(t;x_0))\cdot
g_2(\varphi_0(t;x_0),t+t_0)\\
&\quad-f_2(\varphi_0(t;x_0))\cdot g_1(\varphi_0(t;x_0),t+t_0)\Big)dt.
\end{aligned}
\label{eq:Melnikov}
\end{equation}
Then, if there exists $\bt_0$ such that
\begin{enumerate}
\item $M(\bt_0)=0$
\item $M'(\bt_0)\neq0$,
\end{enumerate}
then, for $\varepsilon>0$ small enough, the perturbed system possesses a
periodic orbit of period $mT$ with initial condition
$x_\varepsilon=x_0+O(\varepsilon)$ at
$t=\bt_0:$
\begin{equation*}
\varphi_\varepsilon(\bt_0+mT;\bt_0,x_\varepsilon)=x_\varepsilon.
\end{equation*}
\end{theorem}
\begin{remark}
The Melnikov function does not require to compute the perturbed
flow. One just needs to evaluated the cross product of the unperturbed
field and the perturbation along a solution of the unperturbed system,
$\varphi_0(t;x_0)$.
\end{remark}
\begin{remark}\label{rem:election_x0}
Given $n$, $m$ and $T_c=\frac{mT}{n}$, the Melnikov function depends
on the election of $x_0$ at the periodic orbit $H(u,v)=c$.
\end{remark}
\begin{remark}
The Melnikov function is $mT$-periodic.
\end{remark}
\begin{exercise}\label{ex:analytical_melnikov}
Compute the analytical expression of the integrand of the Melnikov
function for system~\eqref{eq:ode_pendulum} (you don't have to do the
integral!).
\end{exercise}
\begin{exercise}\label{ex:melnikov}
Write a little program in Matlab to numerically compute the Melnikov
integral obtained in Exercise~\ref{ex:analytical_melnikov}.
\begin{enumerate}
\item Try different $n$'s and $m$'s. Do you always observe zeros of
the Melnikov function?
\item Try a different perturbation with higher harmonics, for example
$g(t)=\sin(\omega t)+4\cos(2\omega t)$. What do you observe?
\end{enumerate}
\end{exercise}
\begin{solution}
The function {\em melnikov.m} performs the Melnikov integral. However,
the initial condition for the unperturbed periodic orbit has to be
properly set. One option is to proceed as in previous exercises: fix a
periodic orbit of the unperturbed system given by some energy level
$H(u,v)=c\in(-1,1)$ by fixing its initial condition at the vertical
axis: $H(0,v_0)=\frac{v_0^2}{2}-1=c$. Then, using the function {\em period.m}
from Exercise~\ref{ex:period} we compute $T_c$. Given $n$ and $m$, the
period $T$ of the forcing becomes $T=nT_c/m$.\\
However, in real applications, we may have a given fixed $T$. In this case,
we need to find the level of energy $c$ corresponding to $T_c=mT/n$.
This could be done by combining the function {\em period.m} from
Exercise~\ref{ex:period} and the Matlab routine {\em fzero} to obtain
the initial condition of an $nT/m$ periodic orbit of the unperturbed
system.\\
Anyhow, once we have $n$, $m$, $T$ and $x_0$ we compute the Melnikov
integral of Theorem~\ref{theo:Melnikov}.
\begin{figure}
\begin{center}
\subfigure[\label{fig:melnk_func1}]{
\includegraphics[width=0.7\textwidth]{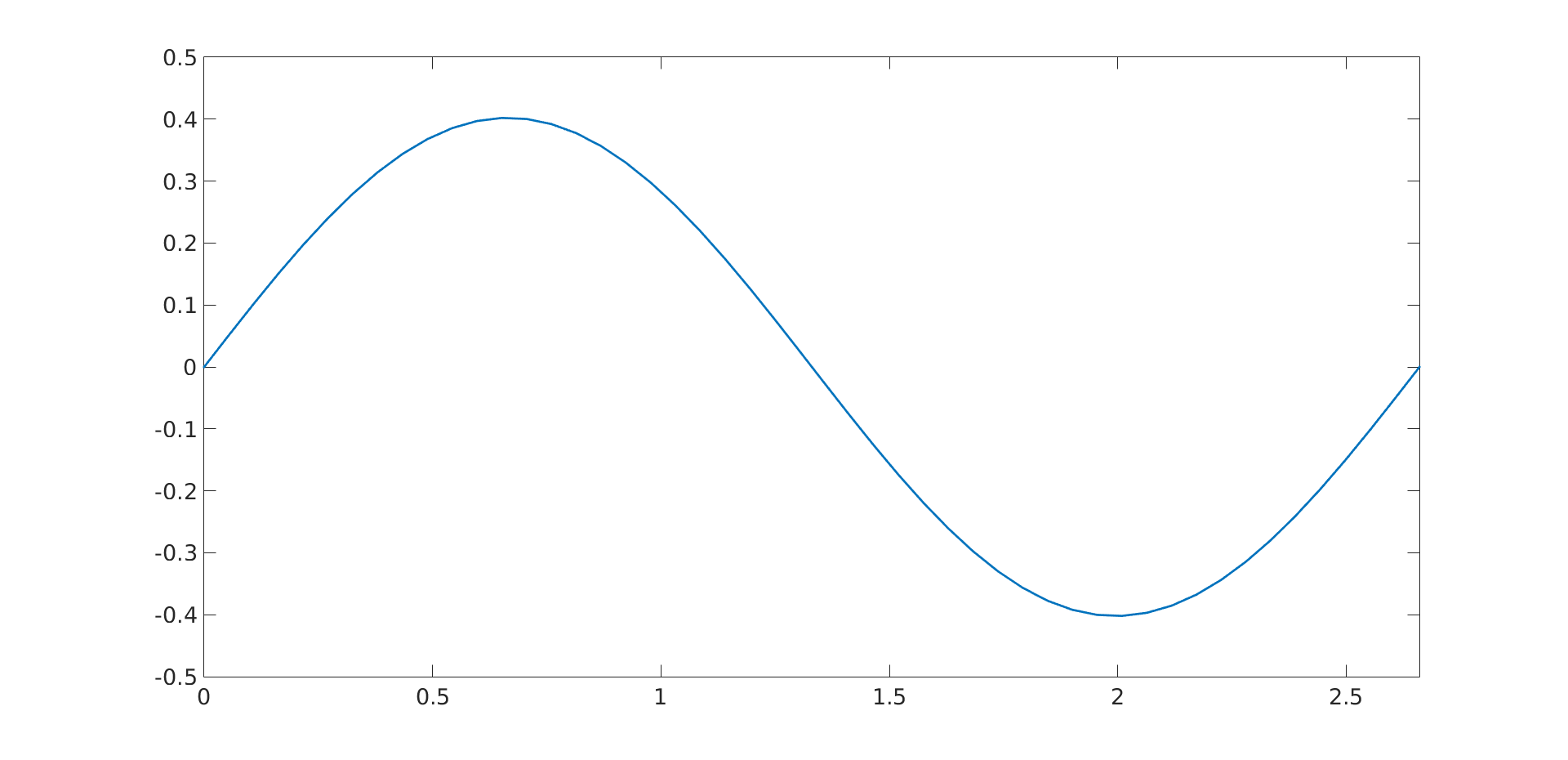}}
\subfigure[\label{fig:melnk_func2}]{
\includegraphics[width=0.7\textwidth]{melnk_func1}}
\subfigure[\label{fig:melnk_func3}]{
\includegraphics[width=0.7\textwidth]{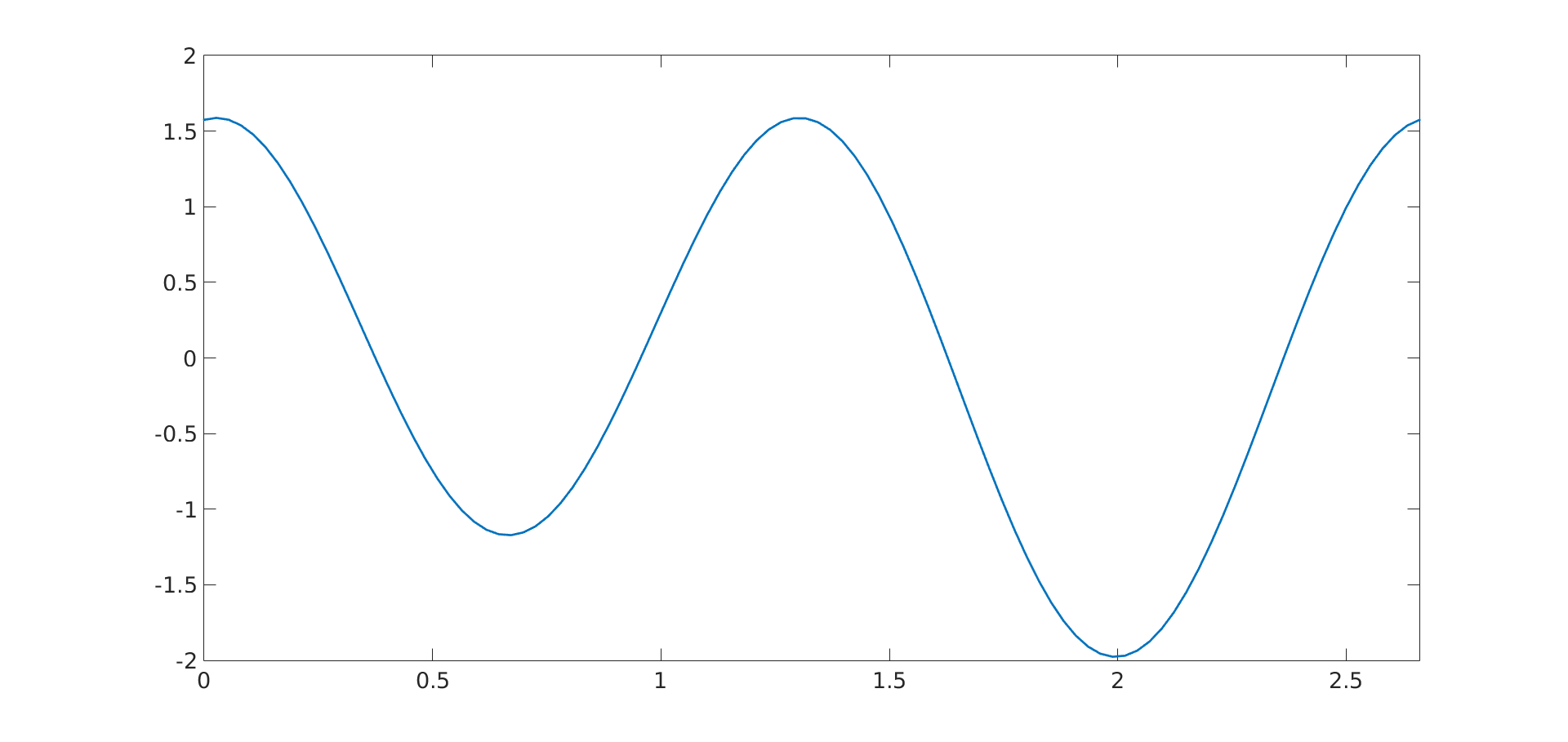}}
\end{center}
\caption{Melnikov function for $g(x,t)=\sin(\omega t)$ (a) and
$g(x,t)=\sin(\omega t)+20\cos(2\omega t)$ (b) and $g(x,t)=\sin(\omega
t)+10^7\cos(2\omega t)$ (c).}
\label{fig:melnikov_functions}
\end{figure}
In Figure~\ref{fig:melnikov_functions} we show the Melnikov function
for three different perturbations for $m=3$ and $n=1$. Note that
Figures~\ref{fig:melnk_func1} and~\ref{fig:melnk_func2} are almost the
same even if the second one contains a higher harmonic. However, as
shown in Figure~\ref{fig:melnk_func3}, if the weight of the second
harmonic is drastically increased, then new zeros appear leading to
more periodic orbits. Note however that, as this weight is so high,
the value of $\varepsilon$ for which one can observe the new periodic
orbits may be extremely small.

When considering the subharmonic case $n=2$ and $m=3$ corresponding to
rotation number $2/3$, it turns out that the Melnikov function is
identically zero for $g(x,t)=\sin(\omega t)$. Therefore, one cannot
prove nor discard the existence of $3$-periodic orbit with rotation
number $2/3$ for such a perturbation and a second order analysis is
required. However, as shown in Figure~\ref{fig:no_2-3_po} when
simulating the stroboscopic map, there is no evidence of existence of
such a periodic orbit.\\
\begin{figure}
\begin{center}
\includegraphics[width=1\textwidth]{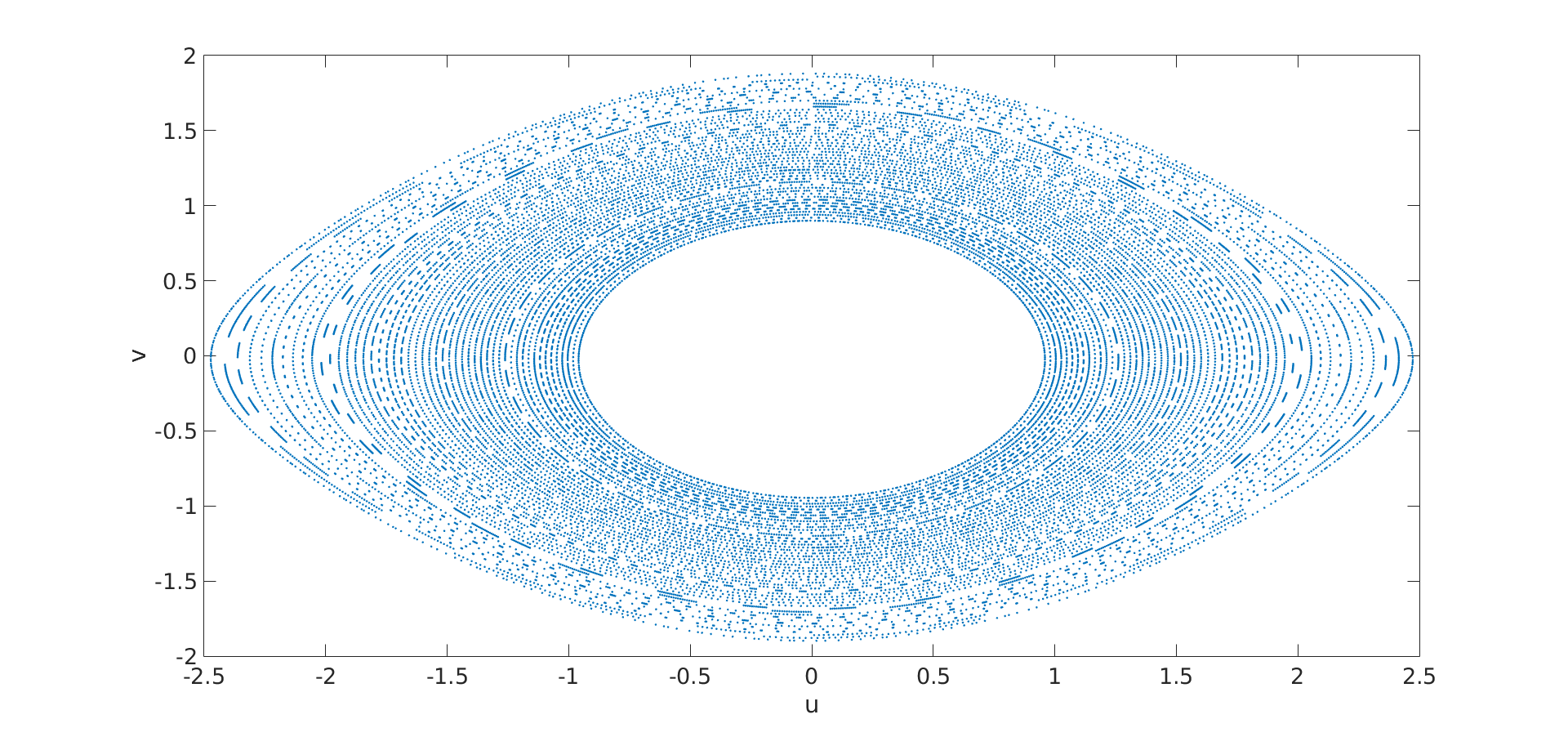}
\end{center}
\caption{Dynamics of the stroboscopic map for $g(x,t)=\sin(\omega t)$
and $\varepsilon=0.01$ around the unpertubed periodic orbit with
rotation number $2/3$. Such periodic orbit does not seem to persist
for $\varepsilon>0$.}
\label{fig:no_2-3_po}
\end{figure}
\begin{figure}
\begin{center}
\includegraphics[width=1\textwidth]{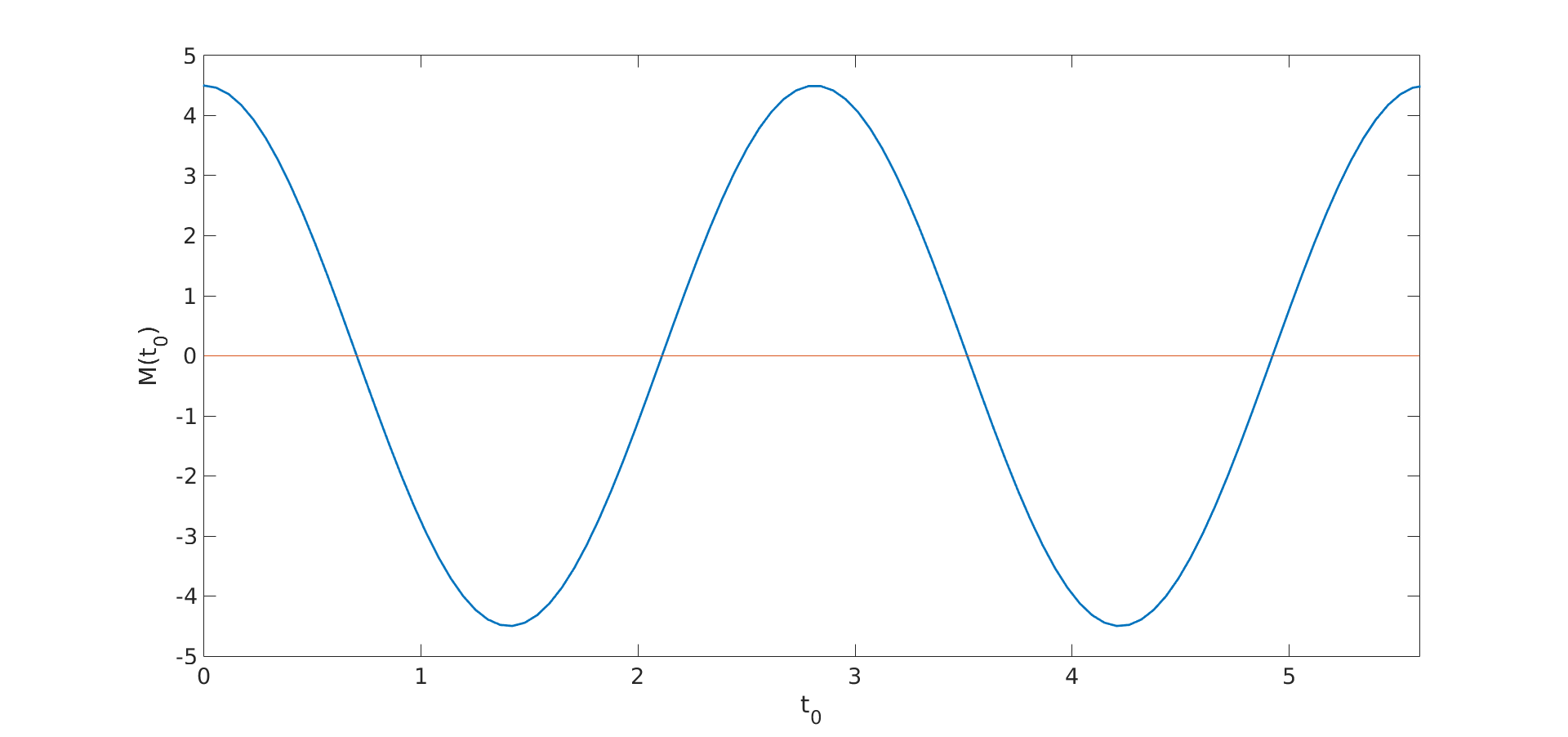}
\end{center}
\caption{Melnikov function for $g(x,t)=\sin(\omega t)+4\cos(2\omega
t)$ using $x_0=(0,1.7)$, $\omega=\frac{2\pi}{T}$, $T=\frac{2T_c}{3}$
and $c=H(0,1.7)$.}
\label{fig:melnkv_2-3}
\end{figure} 
However, as shown in Figure~\ref{fig:melnkv_2-3}, when using a
pertubation with higher harmonics such as
$g(x,t)=\sin(\omega)+4\cos(2\omega t)$, the Melnikov function
possesses simple zeros.
\end{solution}
\subsection{Interpretation}
As we have seen in Section~\ref{sec:unpert_strobo}, when
$\varepsilon=0$, when $T_c=mT/n$, all points satisfying $H(u,v)=c$ are
$m$-periodic points of the stroboscopic map. However, for
$\varepsilon>0$, only some of these points may persist as $m$-periodic
points of the perturbed stroboscopic map. Therefore the following
question arises: given an $m$-periodic point $x_0$ satisfying
$H(x_0)=c$ does this point persist after the perturbation? In other
words, does it exist some point $x_\varepsilon$ $\varepsilon$-close to
$x_0$ such that
\begin{equation*}
s^m(x_\varepsilon)=x_\varepsilon?
\end{equation*}
Taking into account Theorem~\ref{theo:Melnikov}, the answer is that
this depends on $t_0$. That is, once we have chosen $x_0$, the
(simple) zeros of the Melnikov function provides us the proper values
of $t_0$ such that $x_\varepsilon$ exists and is $\varepsilon$-close
to $x_0$.\\
Recalling exercise~\ref{ex:perturbed_strobo_map}, when $\varepsilon>0$
(small), the stroboscopic may exhibit isolated periodic points. These
points are $\varepsilon$-close to the curve given by $H(u,v)=c$
(periodic orbit of the unperturbed system); however, these points
twist when $t_0$ is varied. The zeros of the Melnikov function tell us
how much do we have to twist in order to have one of these points
$\varepsilon$-close to the chosen $x_0$.
\section{Computation of periodic orbits using the stroboscopic map}\label{sec:compute_po_strobo}
We now want to numerically compute the initial condition
($x_\varepsilon$) for a periodic orbit given by
Theorem~\ref{theo:Melnikov}. Here we describe a Newton method to do
that in a general setting: $n$-dimensional not necessary Hamiltonian
system.

Assume we want to compute a periodic orbit of an $n$-dimensional
non-autonomous periodic system of the form
\begin{equation}
\dot{x}=f(x,t),
\label{eq:field}
\end{equation}
where $f(x,t)$ is a $T$-periodic field in $t$: $f(x,t+T)=f(x,t)$, for
any $x\in\RR^n$.\\
Due to the periodicity, instead of using Poincar\'e maps in the state
space it will be much more convenient to use the so-called
stroboscopic map, which is indeed a Poincar\'e map in the extended state
space (adding time as a variable and using a section in time). This
map is given by flowing system~\eqref{eq:field} for a time $T$ with
initial condition $(x_0,t_0)$:
\begin{equation*}
s(x_0)=\varphi(t_0+T;x_0,t_0).
\end{equation*}
Provided that system~\eqref{eq:field} is $T$-periodic in $t$,
$s(x)$ becomes a map from the time section $t=t_0$ to itself:
\begin{equation*}
\begin{array}{cccc}
s:&\Sigma_{t_0}&\longrightarrow &\Sigma_{t_0}\\
&x&\longmapsto&\varphi(t_0+T;x,t_0),
\end{array}
\end{equation*}
where
\begin{equation*}
\Sigma_{t_0}=\left\{ (x,t)\in \RR^2\times \Tt
\right\},\,\Tt:=\RR/T\mathbb{Z}.
\end{equation*}
Imagine that we want to compute the initial condition, $x^p$ for a periodic
orbit of system~\eqref{eq:field} at the section $\Sigma_{t_0}$.
Recalling its periodicity, the period of such a periodic orbit must be
a multiple of $T$, say $mT$. In other words, we must look for periodic
points of the map $s$, that is, points $x^p$ such that
\begin{equation*}
s^m(x^p)=x^p.
\end{equation*}
One of the most extended methods for numerically computing such points
is the Newton method, assuming that one has some idea about where such
a points lies, as  we need a first approximation for the Newton method
to converge where we want.
\subsection{Newton method for fixed points}
Let us assume that we are looking for a $T$-periodic orbit; that is,
$m=1$ and we want to find a fixed point of the map $s$. Then, we want
to solve the equation
\begin{equation*}
F(x):=s(x)-x=0.
\end{equation*}
The Newton method consists of considering the linear approximation of
the function $F$ around some point $(x_0,F(x_0))$ (which is a first
approximation of the solution we are looking for) and solve the linear
system instead. This provides a new point, $(x_1,F(x_1))$ which,
hopefully, is a more accurate solution than the initial one.\\
The linear approximation of equation $F(x)=0$ around $(x_0,F(x_0))$ becomes
\begin{equation}
F(x_0)+DF(x_0)(x-x_0)=0,
\label{eq:approximation1}
\end{equation}
where
\begin{equation*}
DF(x_0)=Ds(x_0)-I.
\end{equation*}
In Section~\ref{sec:variational} we will see how to compute the
differential $D_xs(x)$.\\
If we solve Equation~\eqref{eq:approximation1} for $x$ we get
\begin{equation*}
x=x_0-\big(DF(x_0)\big)^{-1}F(x_0).
\end{equation*}
This leads to an iterative process
\begin{equation}
x_{i+1}=x_{i}-(DF(x_{i}))^{-1}F(x_{i}),
\label{eq:newton}
\end{equation}
which converges quadratically to $x^p$ provided that $x_0$ is a good
enough approximation.\\
Alternatively, some programing environments (like Matlab) offer
routines to solve linear equations which might be more efficient than
computing the inverse $(DF(x_0))^{-1}$. In this case, we can consider
the linear equation
\begin{equation}
DF(x_i)x_{i+1}=DF(x_i)x_i-F(x_i),
\end{equation}
which we want to solve for $x_{i+1}$.\\
\begin{remark}\label{rem:DF_invertibility}
Note that the Newton method requires the matrix $DF(x_i)$ to be
invertible. This implies two things:
\begin{itemize}
\item $DF$ must be invertible at the starting point
\item $DF$ must be invertible at the fixed point! This implies that
the Newton method will have troubles if the fixed point we are looking
for has a eigenvalues real eigenvalue equal to $1$. Note that, as $Ds$
and $Id$ commute, the eigenvalues of $DF$ are the one of $Ds$ minus
one.  In two dimensions, this is degenerate and pathological, but
possible. In dimension one this occurs frequently.
\end{itemize}
\end{remark}
\subsection{Newton method for periodic orbits}\label{sec:differential_po}
Similarly we can apply the Newton method to solve the equation
\begin{equation*}
F(x)=s^m(x)-x
\end{equation*}
to get the same expression as in Eq.~\eqref{eq:newton}. However, in
this case the computation of $DF$ becomes now a bit more tricky. Using
that $s^i(x)=s(s^{i-1}(x))$, we apply the chain rule to get
\begin{equation*}
DF(x_{i})=\prod_{j=i}^{j=1}Ds (s^{j-1}(s)),
\end{equation*}
se we need to evaluate the differential $Ds(x)$ at the points
$s^j(x)$ for $j=1..i$.\\

Although there is nothing wrong with this approach from the
theoretical point of view, in next section we will see a numerical
method to compute $Ds(x)$ which makes the computation of $Ds^m(x)$
straightforward, without needing to multiply matrices (see
Remark~\ref{rem:varia_nperiodic} below).
\subsection{Computation of the differential of the stroboscopic map: the variational equations}\label{sec:variational}
Now the question arises, how do we compute $Ds$? Note that the flow
$\varphi(t;x_0,t_0)$ is straightforward to differentiate with respect
to $t$, as one recovers the field $f$, but we need to differentiate it
with respect to the initial condition $x_0$! But we can do the
following. Applying the fundamental theorem of calculus and the
definition of the flow, we can write
\begin{align*}
s(x_0)&=\varphi(t_0+T;x_0,t_0)=x_0+\int_{t_0}^{t_0+T}\frac{d}{dt}\varphi(t;x_0,t_0)dt\\
&=x_0+\int_{t_0}^{t_0+T}f(\varphi(t;x_0,t_0),t)dt.
\end{align*}
If we now differentiate with respect to $x_0$, we get
\begin{equation}
Ds(x_0)=I_n+\int_{t_0}^{t_0+T}Df\left( \varphi(t;x_0,t_0),t
\right)\cdot D_{x_0}\varphi(t;x_0,t_0)dt,
\label{eq:Ds}
\end{equation}
where $I_n$ is the $n\times n$ identity matrix and we write $D_{x_0}$
to emphasize that we differentiate with respect to $x_0$, while $Df$
is the Jacobian of $f$. Again, by applying the fundamental theorem of
calculus backwards, we realize that Equation~\eqref{eq:Ds} is the
solution of the differential equation
\begin{equation}
\frac{d}{dt}D_{x_0}\varphi(t;x_0,t_0)=Df(\varphi(t;x_0,t_0))\cdot D_{x_0}\varphi(t;x_0,t_0),
\label{eq:variational_equation}
\end{equation}
at $t_0+T$. Equation~\eqref{eq:variational_equation} is called the
(first) variational equation.\\
There is however a more straightforward argument to obtain the
variational equations. If we just differentiate
\begin{equation*}
\frac{d}{dt}D_{x_0}\varphi(t;x_0,t_0),
\end{equation*}
by commuting $Dx$ and $d/dt$ (which is possible of $\varphi$ is at
least $C^2$ with respect to $t$ and $x_0$), we obtain
\begin{equation*}
\frac{d}{dt}D_{x_0}\varphi(t;x_0,t_0)=
D_{x_0}\frac{d}{dt}\varphi(t;x_0,t_0)=D_{x_0}f\left( \varphi(t;x_0,t_0) \right).
\end{equation*}
Applying the chain rule we obtain
\begin{equation*}
\frac{d}{dt}D_{x_0}\varphi(t;x_0,t_0)=
Df\left( \varphi(t;x_0,t_0) \right)\cdot D_{x_0}\varphi(t;x_0,t_0)
\end{equation*}
\noindent Some remarks:
\begin{remark}
If $f(x,t)$ is a field in $\RR^n$, then this equation becomes
and $n\times n$-dimensional differential equation.
\end{remark}
\begin{remark}
Equation~\eqref{eq:variational_equation} is evaluated along the
flow $\varphi(t;x_0,t_0)$, which is unknown. Hence, this equation needs
to be solved together with the equation $\dot{x}=f(x,t)$, leading to a
total system of dimension $n+n\times n$ with initial condition $(x,I_n)$.
\end{remark}
\begin{remark}\label{rem:varia_nperiodic}
If we want to compute $Ds^m(x)$, we just need to integrate the
variational equations from $t_0$ to $t_0+mT$!
\end{remark}

\begin{exercise}
Obtain the variational equations for the perturbed pendulum.
\end{exercise}

\subsection{Combining Newton and Melnikov's method}\label{sec:new-meln}
Let us now consider a system as in Equation~\eqref{eq:field_melnikov}
satisfying conditions~\melnikovconds{}.  As seen in
Sections~\ref{sec:quasi_integrable_dynamics},  some of the periodic
orbits that exist for $\varepsilon=0$ (given as curves $H(u,v)=c$) may
persist after the perturbation. In particular, one needs congruency
between the period of the perturbation and the unperturbed periodic
orbit: $nT_c=mT$ for some $m$ and $n$ co-prime integers.\\
When computing periodic points for $\varepsilon>0$ using the Newton
method described in Section~\ref{sec:differential_po}, it becomes
natural to use some $x_0$ satisfying $H(x_0)=c$ as a seed for the
Newton method, as these periodic points are $\varepsilon$-close to the
curve given by $H(u,v)=c$.  However, unlike in the unperturbed case,
periodic points of the perturbed stroboscopic map become isolated and
hence not any $x_0$ with $H(x_0)=c$ will make the Newton method
converge to the desired point. However, recalling that the system is
non-autonomous, their location depends on the election of $t_0$. This
is where the Melnikov method becomes constructive by providing a
proper seed. More precisely, the Melnikov method purposes to choose
$x_0$ arbitrarily at the curve $H(u,v)=c$ and then use a proper value
of $t_0$ ($\bt_0$ a zero of $M(t_0)$) such that the desired periodic
point, $x_\varepsilon$, is $\varepsilon$-close to $x_0$. Hence, if
$\varepsilon>0$ is small enough and $t_0=\bt_0$ is chosen to be a
simple root of the Melnikov function, the Newton method will
quadratically converge to $x_\varepsilon$ when using $x_0$ as initial
guess.

\begin{figure}
\begin{center}
\subfigure[]{
\includegraphics[width=0.8\textwidth]{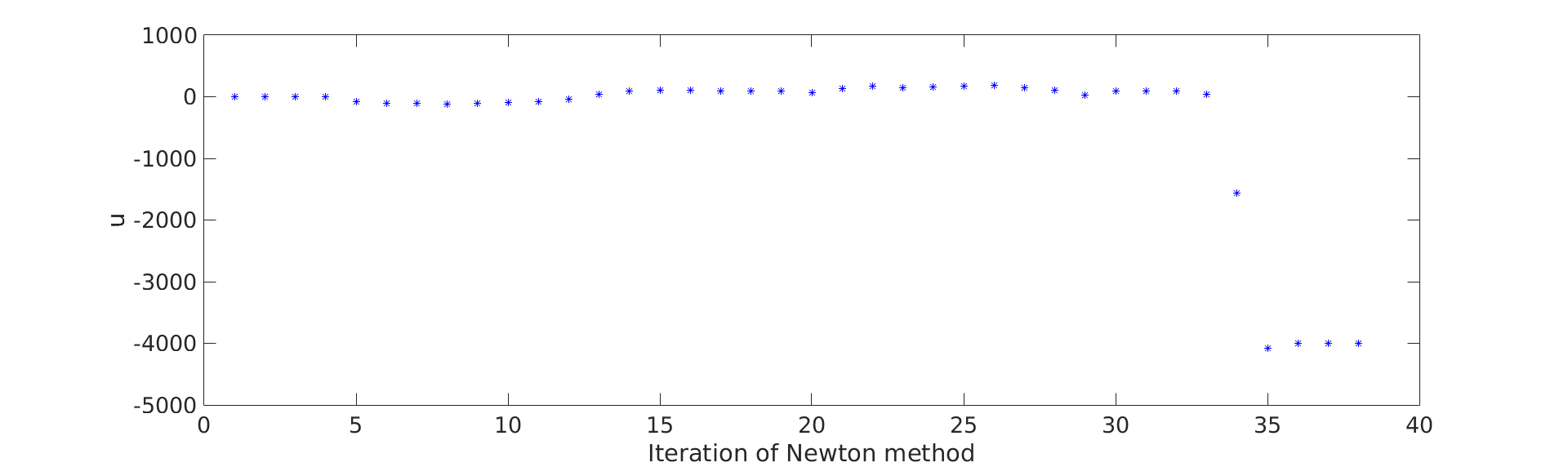}}
\subfigure[]{
\includegraphics[width=0.8\textwidth]{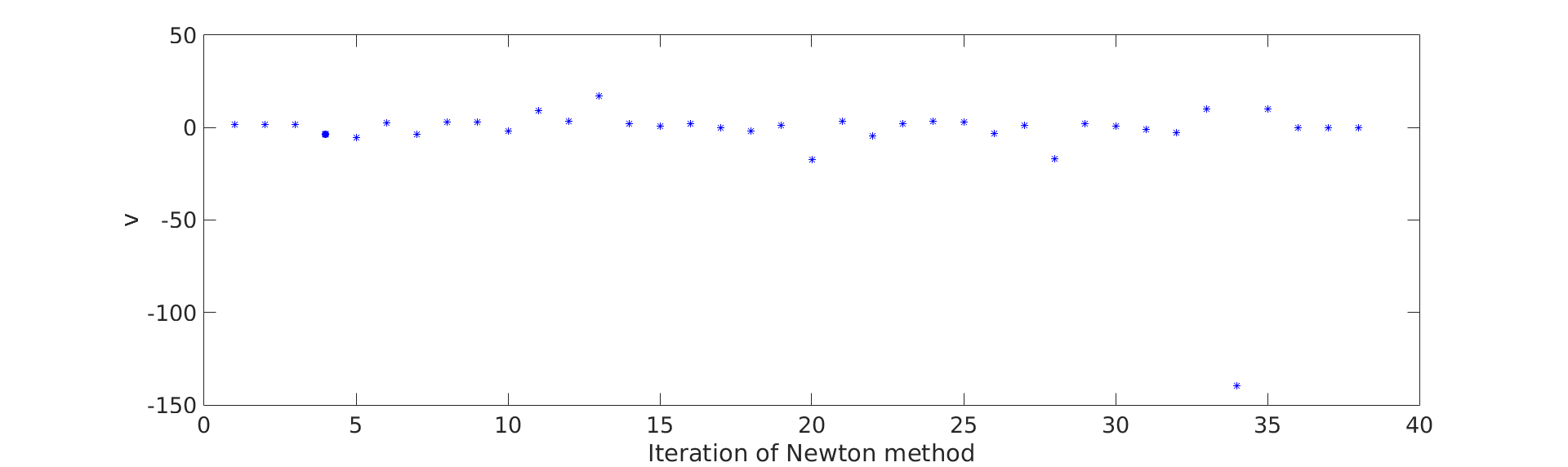}}
\end{center}
\caption{Diverging Newton method for $\varepsilon=0.05$,
initial seed $x_0=(0,1.6)$, $t_0=0$, $T=\frac{T_c}{3}$ and
$c=H(0,1.6)$.}
\label{fig:newton_diverge}
\end{figure}

\begin{figure}
\begin{center}
\subfigure[]{
\includegraphics[width=0.8\textwidth]{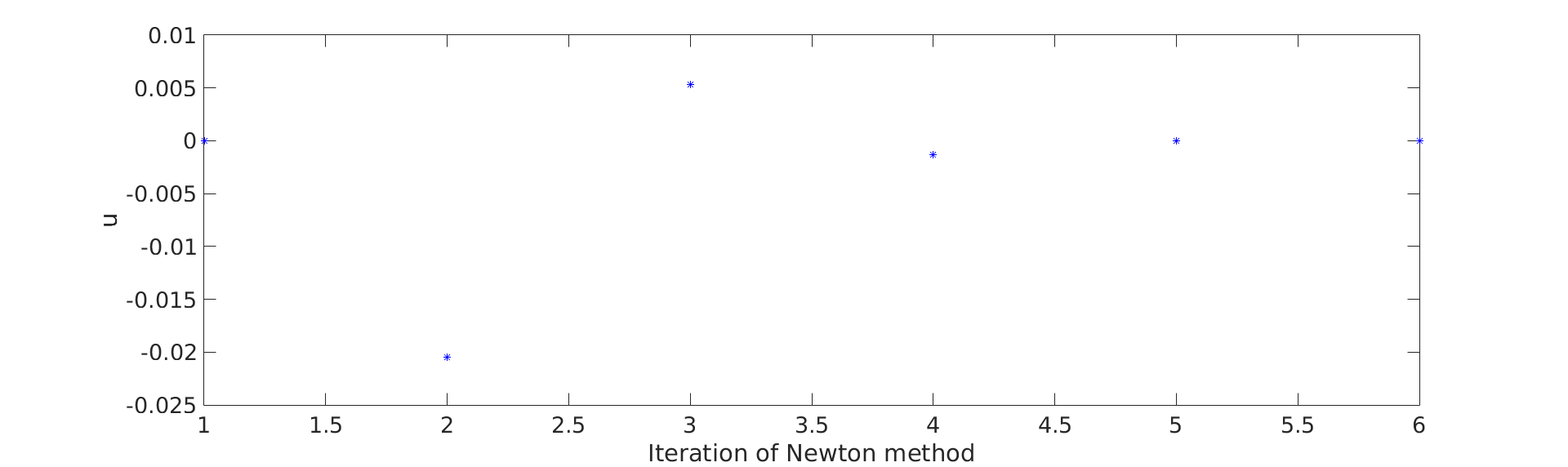}}
\subfigure[]{
\includegraphics[width=0.8\textwidth]{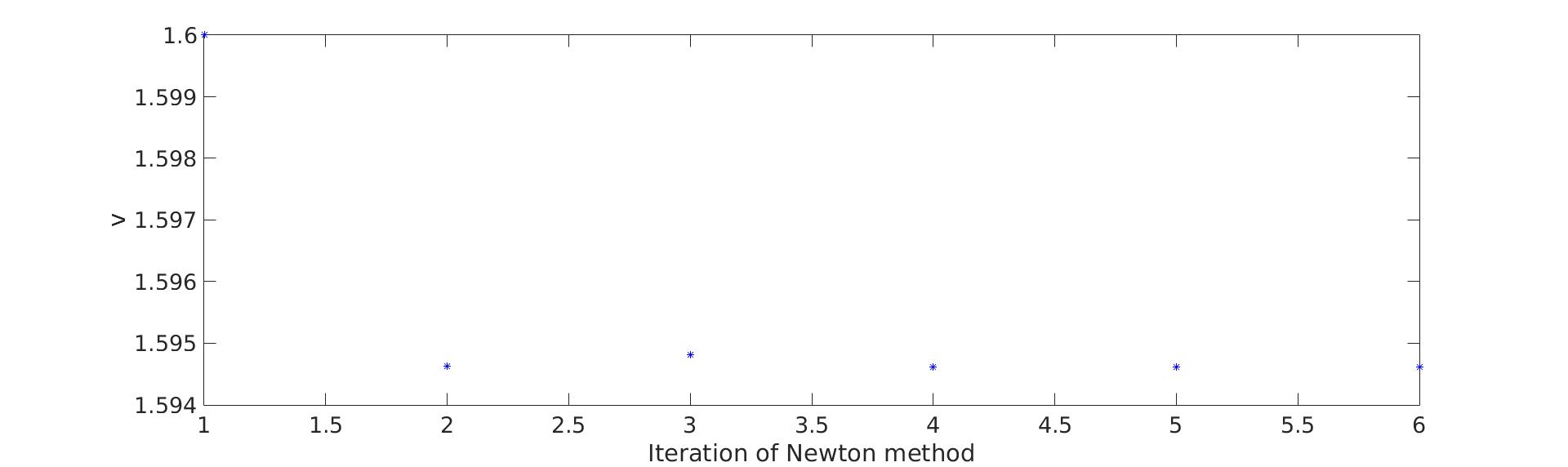}}
\end{center}
\caption{Iterates of the Newton method converging to a saddle periodic
orbit with $\varepsilon=0.01$ and the rest values as in
Figure~\ref{fig:newton_diverge}.}
\label{fig:newton_eps0d01_iterates}
\end{figure}

\begin{figure}
\begin{center}
\includegraphics[width=1\textwidth]{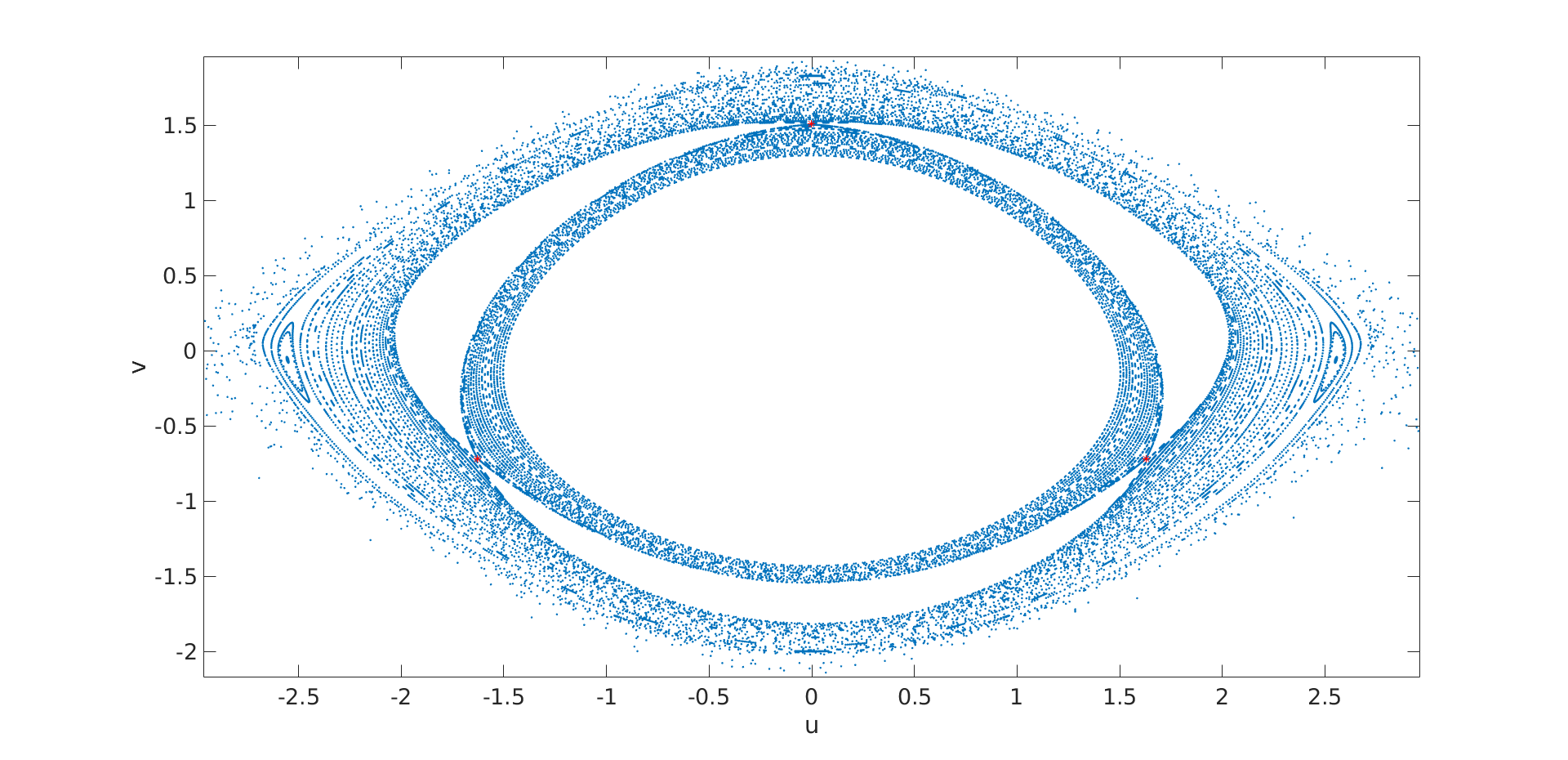}
\end{center}
\caption{Saddle periodic orbit with rotation number $1/3$ for
$\varepsilon=0.17$. The red point close the vertical axis is obtained
by continuation of the Newton method. The other two points ar its
iterates.}
\label{fig:pert3po_saddle}
\end{figure}

\begin{figure}
\begin{center}
\includegraphics[width=1\textwidth]{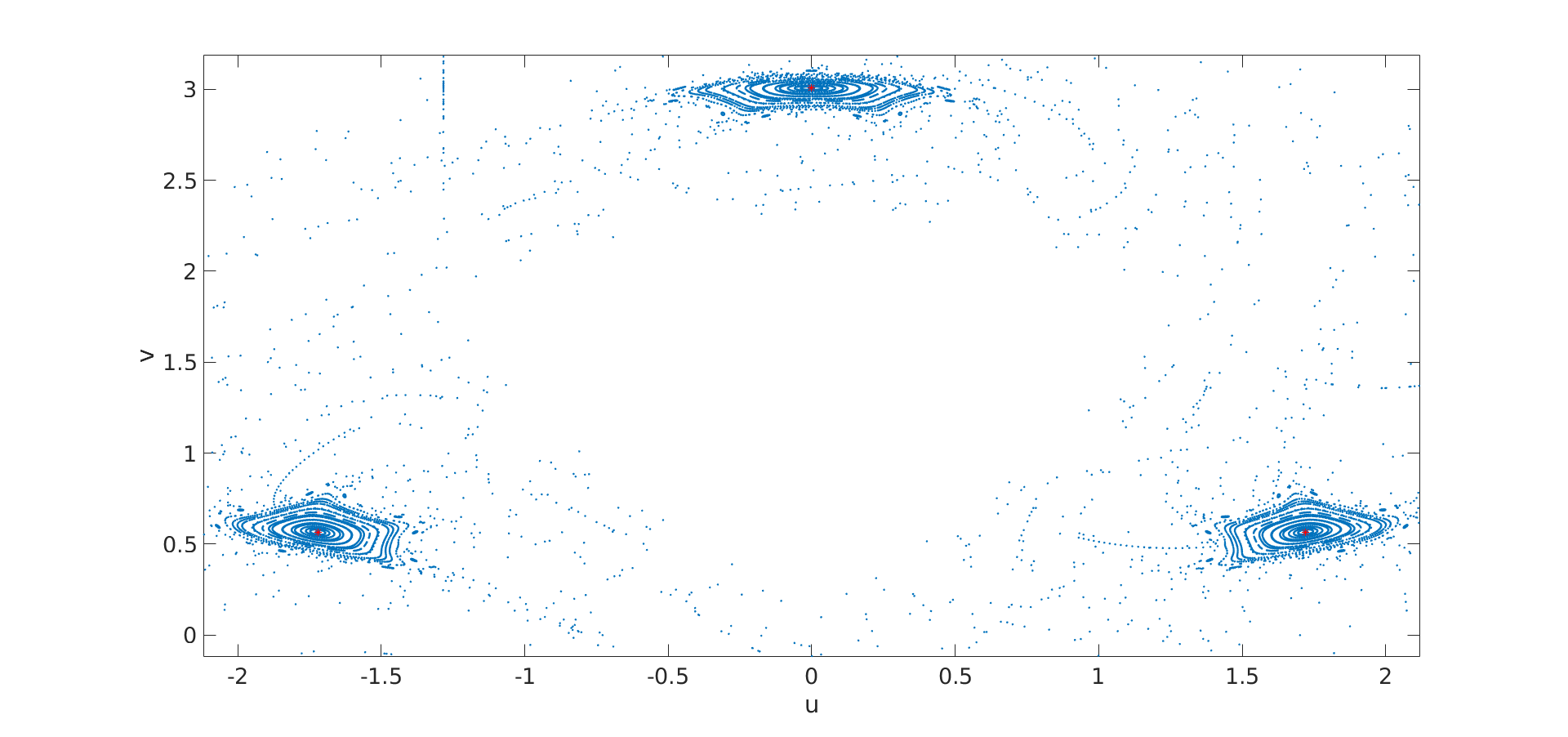}
\end{center}
\caption{Elliptic $3$-periodic orbit for $\varepsilon=2.9$ and
$t_0=\bt_0=T/2\simeq1.33$.}
\label{fig:pert3po_elliptic}
\end{figure}

\begin{figure}
\begin{center}
\includegraphics[width=1\textwidth]{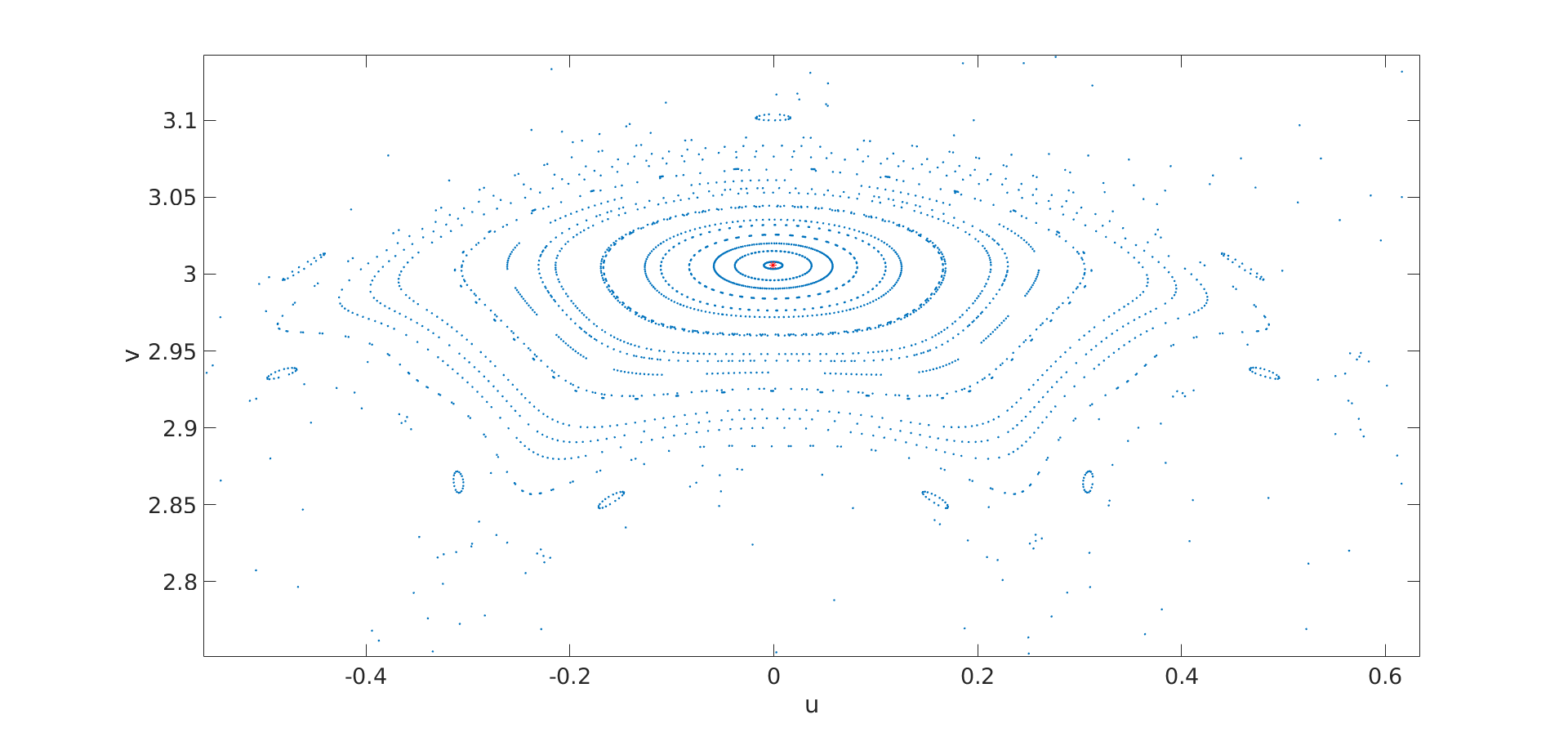}
\end{center}
\caption{Zoom of Figure~\ref{fig:pert3po_elliptic}. In blue the
dynamics of the stroboscopic map, in red the iterates of the periodic
orbit found using the Newton method.}
\label{fig:pert3po_elliptic_zoom}
\end{figure}

\begin{figure}
\begin{center}
\includegraphics[width=1\textwidth]{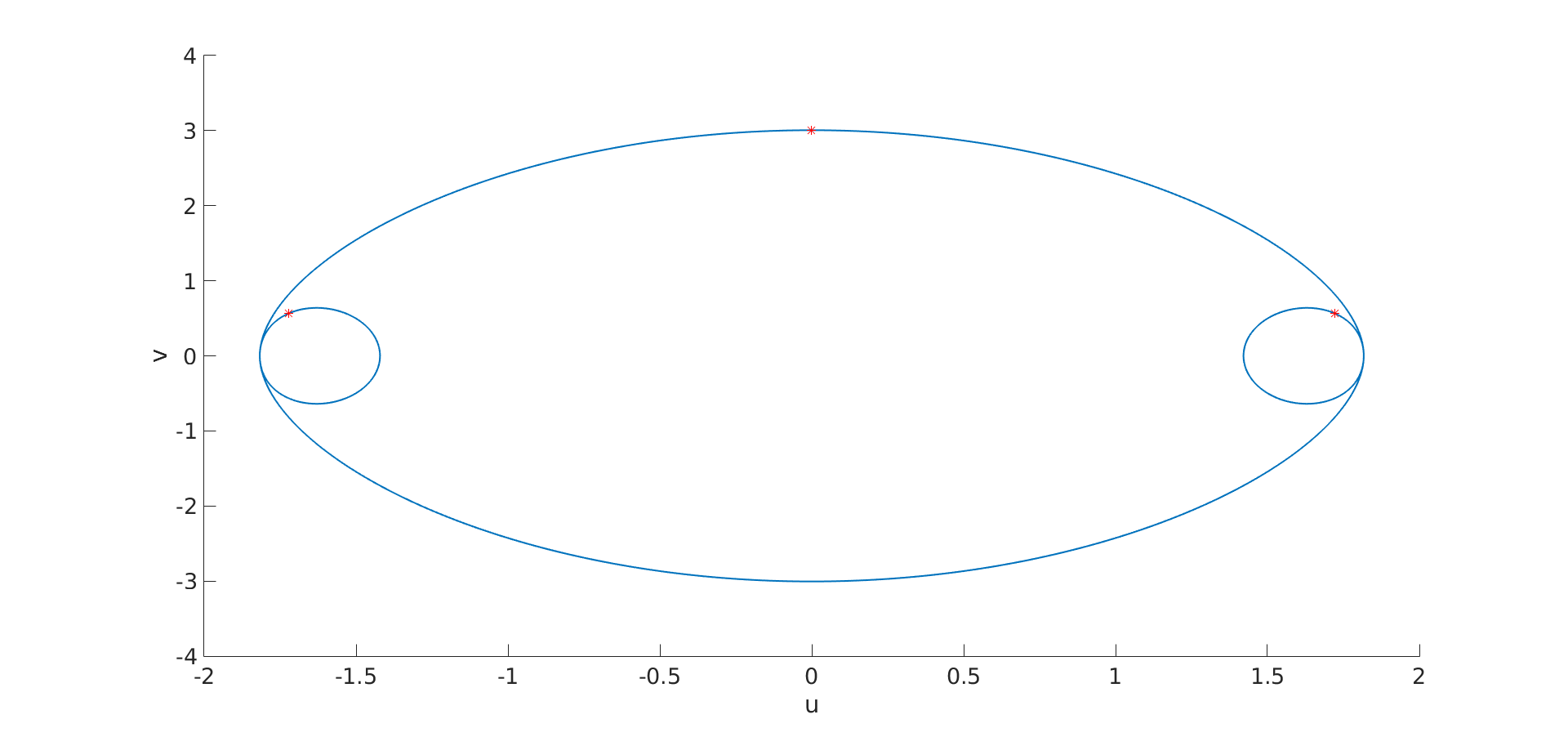}
\end{center}
\caption{Time-continuous trajectory for $t\in[T/2,T/2+3T]$
corresponding to the $3$-periodic orbit shown in
Figure~\ref{fig:pert3po_elliptic} (blue curve). Red points are the
iterates of the stroboscopic map.}
\label{fig:pert3po_elliptic_TC}
\end{figure}

\begin{figure}
\begin{center}
\includegraphics[width=1\textwidth]{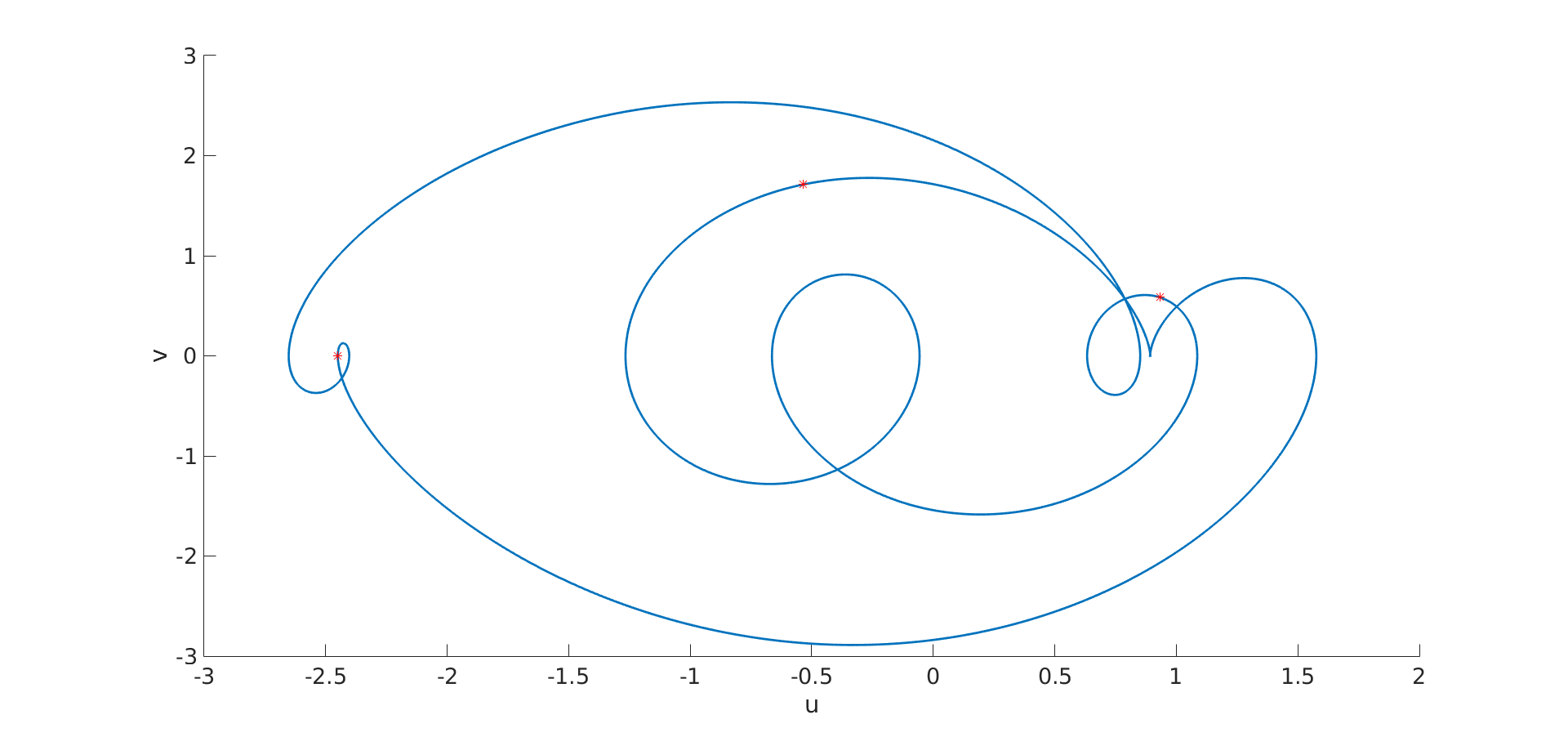}
\end{center}
\caption{Periodic orbit with rotation number $2/3$ for
$g(x,t)=\sin(\omega t)+4\cos(2\omega t)$ and $\varepsilon=0.5$.}
\label{fig:pert3po_2-3_TC}
\end{figure}

\begin{exercise}
Write a program in Matlab in order to compute $nT/m$-periodic orbits
for $\varepsilon>0$.
\begin{enumerate}
\item Obtain the variational equations and write them in the file
{\em pertpend\_vari.m}.
\item Implement the Newton method to obtain a fixed point of the
stroboscopoic map. You will need to tune $\omega$ to guarantee that it
exists. In the file {\em Newton.m} you have some lines that you can
use. Note that the initial seed will be taken from
a zero of the Melnikov function, which may have several. First use an
initial seed obtained manually to make sure that your program works.
\item Perform a veeeery little modification in your program to compute
$nT$-periodic orbits instead of a $T$-periodic orbit.  Compute the
eigenvalues of $Ds$ at these periodic orbits. What is their type?
\end{enumerate}
\end{exercise}
\begin{solution}
The variational equations are easy to obtain.\\
The routine {\em Newton.m} contains one step of the Newton method. By
repeating this step until the norm of $F$ is small enough, we we
perform a Newton method for a fixed point. Just by replacing $T$ by
$mT$, the program will compute an $m$-periodic orbit instead of a
fixed point.\\
Now we use the program. We choose to compute a $3$-periodic orbit with
rotation number $1/3$. Hence, we use as seed a point satisfying
$H(x_0)=c$, form some $c$, and use $T=\frac{T_c}{3}$ and
$\omega=\frac{2\pi}{T}$. We choose for
example value we used in previous exercises: $x_0=(0,1.6)$ and
$c$ becomes $c=\frac{1.6^2}{2}-1$. The Melnikov function for
$g(x,t)=\sin(\omega t)$ is shown in Figure~\ref{fig:melnk_func1}, and
possess two simple zeros: $\bt_0=0$  and $\bt_0=T/2\simeq 1.4$.
Using $x_0=(0,1.6)$ and $t_0=0$ as seed for the Newton method, we show
in Figure~\ref{fig:newton_diverge} how the method diverges for
$\varepsilon=0.05$. By reducing $\varepsilon$ to $\varepsilon=0.4$ the
Newton method converges to another periodic orbit. By further
decreasing $\varepsilon$ to $\varepsilon=0.01$, the Newton method
finally converges to the desired periodic orbit (a saddle). In
Figure~\ref{fig:newton_eps0d01_iterates} we show how the method converges in 6
iterations reaching an accuracy of $10^{-8}$. By using the obtained
result as seed for a next Newton method we can follow up this periodic
orbit by increasing $\varepsilon$. In Figure~\ref{fig:pert3po_saddle}
we show what we obtain for $\varepsilon=0.17$.

We now focus on the other periodic orbit given by the zero of the
Melnikov function around $\bt_0=T/2\simeq1.33$. Proceeding similarly,
by continuation of the Newton method, we can follow this periodic
orbit (of the elliptic type) up values higher than $\varepsilon=5$. In
Figure~\ref{fig:pert3po_elliptic} we show the iterates of the
stroboscopic map of this periodic orbit for $\varepsilon=2.9$ and the
dynamics surrounding it. As we can see there, very little of the rest
of the dynamics persist for such large value of $\varepsilon$.
However, when zooming in (see Figure~\ref{fig:pert3po_elliptic_zoom}),
we observe secondary tori and secondary periodic orbits around the
elliptic periodic orbit. In Figure~\ref{fig:pert3po_elliptic_TC} we
show the time-continuous evolution of $u$ and $v$ along this periodic
orbit. As can see, the periodic orbit takes $3$ times $T$ to perform
one loop around the origin, hence reflecting the rotation number
$1/3$. 

We finally discuss the case of rotation number $2/3$. As explained in
Exercise~\ref{ex:melnikov}, we need to include higher harmonics in
order to guarantee the existence of such periodic orbits; otherwise,
the Melnikov function becomes identically zero. Proceeding as before,
using $g(x,t)=\sin(\omega t)+3\cos(2\omega t)$,
$\omega=\frac{2\pi}{T}$, $T=\frac{2T_c}{3}$, $c=H(0,1.7)$, we use
$\bt_0=0.7035$ (first zero of Figure~\ref{fig:melnkv_2-3}) and
$x_0=(0,1.7)$ as seed for the Newton method. Starting with a small
$\varepsilon$ and performing continuation, we obtain the periodic
orbit with rotation number $2/3$ shown in
Figure~\ref{fig:pert3po_2-3_TC} for $\varepsilon=0.5$. As can see
there, the periodic orbit performs $2$ loops when closing at
$t=\bt_0+3T$.
\end{solution}

\section{Computation of  periodic orbits using a Poinca\-r\'e map}
As noted in Remark~\ref{rem:DF_invertibility}, the Newton method to
find fixed points (or periodic orbits) of the stroboscopic map may
fail if that one has a real eigenvalue equal to $1$. Alternatively, we can
use a Poincar\'e map using a section in the state space. This method is
more robust in that sense, but, as we will show below, the computation
of the differential becomes slightly more tricky.

Let us consider the Poincar\'e map
\begin{equation}
\begin{array}{cccc}
P_\varepsilon(v,t):&\left\{ u=0 \right\}\times\RR&\longrightarrow&\left\{ u=0
\right\}\times\RR\\
&(0,v_0,t_0)&\longmapsto&(0,\varphi_1(t^*;0,v_0,t_0),t^*),
\end{array}
\label{eq:poinc_map}
\end{equation}
where
\begin{equation}
\varphi_\varepsilon(t;x_0,t_0)=
\left( 
\begin{array}[]{c}
\varphi_1(t;x_0,t_0)\\\varphi_2(t;x_0,t_0)
\end{array}
 \right)
\label{eq:perturbed_flow_notation}
\end{equation}
is the solution of the perturbed system~\eqref{eq:perturbed_pendulum}
with initial condition $x_0$ at $t=t_0$ and $t^*$ is such
that\footnote{Here we assume that we know how to compute $t^*$. It can
be easily done using a Newton method to solve
Equation~\eqref{eq:return_condition}. Below we show how to compute the
necessary derivatives for the Newton method. Otherwise, it can be
computed using the ``Events'' functionality of Matlab.}
\begin{equation}
\varphi_1(t^*;0,v_0,t_0)=0
\label{eq:return_condition}
\end{equation}
for the second time, as we consider a return map in the direction that
the section is left.\\
Recall that system~\eqref{eq:perturbed_pendulum} is non-autonomous.
Hence, the initial time $t_0$ matters and we have to carry it on.
From now on, we will abuse notation in Equation~\eqref{eq:poinc_map}
and omit writing the first coordinate $u$, as it takes the $0$.

Taking into account the periodicity of the non-autonomous
system~\eqref{eq:perturbed_pendulum}, an initial condition at
$(0,v_0)\in\left\{ u=0 \right\}$ 
 for $t=t_0$ of a
periodic orbit of the full system will need to satisfy
\begin{equation}
P^n(v_0,t_0)=\left( 
\begin{array}[]{c}
v_0\\t_0+mT
\end{array}
 \right)
\label{eq:periodic_orbit_poinc-map}
\end{equation}
for some $n,m>0$. Indeed, the integers $n$ and $m$ are the same as in
previous sections, made explicit through the congruency
equation~\eqref{eq:congruency}. Here one clearly see that the role of
$n$ is to count the ``loops'' that a periodic orbit of period
$mT$ makes around the origin.

Applying the Newton method to the equation
\begin{equation}
F(v_0,t_0)=P^n(v_0,t_0)-\left(
\begin{array}[]{c}
v_0\\t_0+mT
\end{array}
 \right)=\left( 
\begin{array}[]{c}
0\\0
\end{array}
 \right),
\label{eq:newton_equation_poinc-map}
\end{equation}
and arguing as in Section~\ref{sec:compute_po_strobo} we get the
iterative process
\begin{equation}
(v_{i+1},t_{i+1})=(v_i,t_i)-(DF(v_i,t_i))^{-1}F(v_i,t_i),
\label{eq:newton_poinc_map}
\end{equation}
where now
\begin{equation*}
DF(v_i,t_i)=DP^n(v_i,t_i)-\left( 
\begin{array}[]{cc}
1&0\\0&1
\end{array}
 \right).
\end{equation*}
We now wonder, how do we compute $DP^n$? Let us see it for $n=1$
first.\\

\noindent Recalling that $\varphi_\varepsilon=\left( \varphi_1,\varphi_2
\right)$, we get
\begin{equation*}
DP(v_0,t_0)=\left( 
\begin{array}[]{cc}
D_{v_0}\varphi_2(t^*;0,v_0,t_0)&
D_{t_0}\varphi_2(t^*;0,v_0,t_0)\\
\frac{\partial}{\partial v_0}t^*&
\frac{\partial}{\partial t_0}t^*\\
\end{array}
 \right).
\end{equation*}
Note that, in the first row, we have written the total derivatives
$D_{v_0}$ and $D_{t_0}$ instead of partial ones, $\partial/\partial
v_0$ and $\partial/\partial_{t_0}$, because $t^*$ actually depends on
$v_0$ and $t_0$ through Equation~\eqref{eq:return_condition}. Let us
compute such total derivatives.\\
For the first one we get
\begin{equation*}
D_{v_0}\varphi_2(t^*;0,v_0,t_0)=\frac{\partial }{\partial
v_{0}}\varphi_2(t^*;0,v_0,t_0)+\varphi_2'(t^*;0,v_0,t_0)\frac{\partial
t^*}{\partial v_0}.
\end{equation*}
Let us now see how to compute all the terms appearing in this equation.\\
On one hand,
\begin{equation*}
\varphi_2'(t^*;0,v_0,t_0)=f_2(\varphi_\varepsilon(t^*;0,v_0,t_0))+\varepsilon
g_2\left(\varphi_\varepsilon
\left(t^*;0,v_0,t_0),t^*\right)\right),
\end{equation*}
where $f_2+\varepsilon g_2$ is the second coordinate of the field
evaluated at the image of the Poincar\'e map.\\
On the other one, $\frac{\partial \varphi_2(t^*;0,v_0,t_0)}{\partial
v_0}$ is given by integrating the variational equations from $t_0$ to
$t^*$, as we are differentiating with respect an initial condition.\\
What about $\frac{\partial t^*}{\partial v_0}$? Recall that, for given
$v_0$ and $t_0$, $t^*$ (which we somehow know how to compute) solves
Equation~\eqref{eq:return_condition}. Therefore, assuming that
\begin{equation}
\varphi_1'(t^*;0,v_0,t_0)\ne0,
\label{eq:IFT_condition}
\end{equation}
we can use the Implicit Function Theorem to get
\begin{align}
\frac{\partial t^*}{\partial
t_0}&=-\frac{\frac{\partial}{\partial
t_0}\varphi_1(t^*;0,v_0,t_0)}{\varphi_1'(t^*;0,v_0,t_0)}\nonumber\\
&=-\frac{\frac{\partial}{\partial
t_0}\varphi_1(t^*;0,v_0,t_0)}{f_1(\varphi_\varepsilon(t^*;0,v_0,t_0))+\varepsilon
g_1(\varphi_\varepsilon(t^*;0,v_0,t_0),t^*)}.
\label{eq:pt/pt0}
\end{align}
Note that, condition~\eqref{eq:IFT_condition} is satisfied, as the
flow is transversal to the section $\left\{ u=0 \right\}$. If it were
tangent, then we would have a problem, of course!\\
The denominator  of the last equation is just the first coordinate of
the perturbed field evaluated at the image of the Poincar\'e map.\\
Regarding $\frac{\partial}{\partial t_0}\varphi_1(t^*;0,v_0,t_0)$,
note that it is a derivative with respect to the initial time, $t_0$.
Adding time as a variable, this can be computed integrating the
variational equations of the system
\begin{equation}
\begin{aligned}
\dot{u}&=v\\
\dot{v}&=-\sin(u)+\varepsilon \sin(\omega s)\\
\dot{s}&=1,
\end{aligned}
\label{eq:extended_system}
\end{equation}
where now $s$ plays the role of time and $t_0$ becomes $s_0$.\\

\noindent Proceeding similarly, the elements of the second row of
$DP$ become
\begin{equation*}
\frac{\partial t^*}{v_0}=-\frac{\frac{\partial}{\partial
v_0}\varphi_1(t^*;0,v_0,t_0)}{f_1(\varphi_\varepsilon(t^*;0,v_0,t_0))+\varepsilon
g_1\left( \varphi_\varepsilon(t^*;0,v_0,t_0) \right)},
\end{equation*}
and $\frac{\partial t^*}{\partial t_0}$ is already given in
Equation~\eqref{eq:pt/pt0}\\

\noindent The differential $DP^n(v_0,t_0)$ can be computed proceeding
similarly as in Section~\ref{sec:differential_po}, mutliplying $DP$
evaluated at the iterates $P^n(v_0,t_0)$ or integrating the previous
variational equations until the $2n$-th impact with the section
$\left\{ u=0 \right\}$ occurs.

\begin{exercise}
Compute analytically the variational equations of
system~\eqref{eq:extended_system}.
\end{exercise}

\begin{exercise}
Write a program in Matlab to compute initial conditions for periodic
orbits by performing a Newton method to
Equation~\eqref{eq:newton_equation_poinc-map}. For that you will need
to use all the artillery you developed in the previous exercises.
\end{exercise}
\def\zh{Zh}\def\yu{Yu}\def\ya{Ya}

\end{document}